\documentclass[preprint,11pt]{elsarticle}
\usepackage[letterpaper,left=1.91cm,right=1.91cm,top=1.91cm,bottom=2cm,]{geometry}
\usepackage[pdftex,colorlinks,linkcolor={blue},citecolor={blue}]{hyperref}

\usepackage{amsmath,amsfonts}
\usepackage{amsthm}
\usepackage{array}
\usepackage[caption=false,font=footnotesize,labelfont=normalfont,textfont=normalfont]{subfig}
\usepackage{url}

\usepackage{graphicx}
\usepackage{cases}
\usepackage{xcolor}
\usepackage{pgfplots}
\pgfplotsset{compat=newest}
\usetikzlibrary{plotmarks}
\usepackage{standalone}
\usepackage{footnote}
\usepackage{multirow}
\usepackage{booktabs}
\usepackage{arydshln}

\usepackage[ruled,vlined,linesnumbered]{algorithm2e}
\SetArgSty{textnormal}

\makeatletter
\def\@author#1{\g@addto@macro\elsauthors{\normalsize%
    \def\baselinestretch{1}%
    \upshape\authorsep#1\unskip\textsuperscript{%
      \ifx\@fnmark\@empty\else\unskip\sep\@fnmark\let\sep=,\fi
      \ifx\@corref\@empty\else\unskip\sep\@corref\let\sep=,\fi
      }%
    \def\authorsep{\unskip,\space}%
    \global\let\@fnmark\@empty
    \global\let\@corref\@empty  
    \global\let\sep\@empty}%
    \@eadauthor={#1}
}
\makeatother

\usepackage{accents}
\newlength{\dhatheight}

\newcommand{\m}[1]{\mathcal{#1}}

\renewcommand{\t}[1]{\text{#1}}

\newtheorem{proposition}{Proposition}
\newtheorem{remark}{Remark}

\theoremstyle{definition}
\newtheorem{definition}{Definition}[section]

\begin{document}

\begin{frontmatter}
\title{Towards a Multimodal Charging Network: Joint Planning of Charging Stations and Battery Swapping Stations for Electrified Ride-Hailing Fleets}

\author{Zhijie Lai}
\ead{zlaiaa@connect.ust.hk}
\author{Sen Li\corref{cor}}
\cortext[cor]{Corresponding author}
\ead{cesli@ust.hk}
\address{Department of Civil and Environmental Engineering, The Hong Kong University of Science and Technology, \\Clear Water Bay, Kowloon, Hong Kong}

\begin{abstract}
This paper considers a multimodal charging network in which charging stations and battery swapping stations are jointly built to support an electric ride-hailing fleet synergistically. Our argument is based on the observation that charging an EV is a time-consuming burden, and battery swapping faces scaling issues due to its deployment costs. However, charging stations are cost-effective, making them ideal for scaling up EV fleets, while battery swapping stations offer quick turnaround and can be deployed in tandem with charging stations to improve fleet utilization and reduce operational costs. To fulfill this vision, we consider a ride-hailing platform that jointly builds charging and battery swapping stations to support an EV fleet. An optimization model is proposed to capture the platform's planning and operational decisions. In particular, the model incorporates essential components such as elastic passenger demand, spatial charging equilibrium, charging and swapping congestion, etc. The overall problem is formulated as a nonconcave program. Instead of pursuing the globally optimal solution, we establish a tight upper bound through relaxation and decomposition, allowing us to evaluate the solution optimality even in the absence of concavity. Through case studies for Manhattan, New York City,  we find that joint planning of charging and battery swapping stations outperforms deploying only one of them, yielding a total profit that is 11.7\% higher than swapping-only deployment under a limited budget, and 17.5\% higher than charging-only deployment under a sufficient budget. These results underscore the complementary benefit between charging and battery swapping facilities.
\end{abstract}

\end{frontmatter}

\section{Introduction}
The electrification of transportation is essential for establishing a green and sustainable mobility future. In recent years, governments worldwide have committed to supporting the adoption of electric vehicles (EVs) through a suite of regulations and subsidies. Among all types of vehicles, commercial fleets (e.g., ride-hailing cars) are viewed as a promising sector for transport electrification due to their higher annual mileage. Companies like Uber and Didi, with their extensive global operations, have the potential to drive significant progress in pollution reduction. Uber alone provides approximately 21 million rides per day across over 10,000 cities \cite{uber_report_2022}, while Didi's global daily trips can even exceed 50 million \cite{didi_report_2022}. However, the vast majority of these rides are delivered by fossil fuel vehicles, with EVs accounting for less than 5\% of Uber drivers' miles driven in Europe \cite{guardian_electric_2021}, and less than 0.5\% in the US \cite{fleming_policy_2020}. In fact, the environmental benefits of electrifying ride-hailing fleets are estimated to be three times higher than those of electrifying private cars \cite{jenn_emissions_2020}. To unlock this enormous potential, California regulators have passed a legislation requiring transportation network companies to transition from petroleum to electric vehicles by 2030 \cite{forbes_california_2021}. In response, both Uber and Lyft have pledged to offer zero-emission mobility-on-demand services by the end of the decade  \cite{uber_together}\cite{lyft_our}.

However, charging infrastructure remains a persistent barrier that hinders the electrification of ride-hailing fleets in the long term. Unlike private cars, commercial fleets need to run on the road as long as possible to accrue profits, while most ride-hailing EV drivers need to top up the batteries once or twice every day due to the limited battery range \cite{sanguinetti_characteristics_2021}. The frequent need for charging, coupled with the longer charging time compared to refueling a gas-powered vehicle, results in significant downtime, incurring a substantial opportunity cost for the drivers \cite{liu_regulatory_2022}. This makes EVs less appealing for ride-hailing use. As per \cite{roni_optimal_2019}, charging time comprises 72-75\% of the total downtime spent on charging trips. This indicates that charging time, rather than travel time or waiting time, is the dominant bottleneck for EV charging, even with a dense charging station network. While fast charging can alleviate this issue, it typically relies on high voltage and draws substantial power from the grid. Massively deploying fast chargers in densely populated urban areas is impractical due to the potential voltage drop and instability it can cause in the already strained power distribution network. Additionally, upgrading the power distribution infrastructure to support such high-power charging can be prohibitively expensive. As a result, plug-in charging stations presents an inherent bottleneck for the ride-hailing industry: regardless of the density of the charging network, there is always a significant downtime for the commercial fleet due to the limited charging speed that can be maximally supported by the capacity of existing power distribution networks in large metropolitan areas. This leads to low fleet utilization and high operational costs for the ride-hailing platform, thereby limiting the economic viability of electrifying ride-hailing fleets.

Battery swapping is an alternative method for recharging depleted EVs, where an empty battery is replaced with a fully charged one in a few minutes. There was a heated debate over plug-in charging and battery swapping, after which many western countries predominantly invested in charging stations instead of battery swapping stations. One primary reason for this preference is the lack of standardized batteries among EV manufacturers. In contrast, China has taken a more impartial stance, encouraging both plug-in charging and battery swapping technologies. This policy is based on the belief that, standardizing batteries is more feasible for commercial vehicles than for privately-owned EVs. In recent years, governments are pushing electrification among taxi fleets, which naturally come with a unified battery model (see Shenzhen \cite{EV_shenzhen} for an example). In the private sector, ride-hailing platforms have actively collaborated with automobile manufacturers to produce affordable EVs tailored for ride-hailing services. Notable examples include the battery-swappable Cao Cao 60 \cite{geely}, a joint effort between Caocao Mobility and Geely, and the D1 \cite{didi}, co-developed by Didi and BYD. Uber is also working with automakers to develop custom-made EVs for ride-hailing and delivery applications \cite{wsj}. Large-scale deployment of these customized EVs can naturally address the compatibility issue, making battery swapping a viable solution. Furthermore, future mobility-on-demand systems are envisioned to be driven by electric and autonomous vehicles that operate under the centralized control of the platform. In fact, both Uber and Didi have announced their plans to offer services using robotaxis \cite{uber_waymo}\cite{didi_robotaxi}. This important trend further justifies a unified battery standard and the feasibility of battery swapping. The advantages of battery swapping can be tremendous for electric ride-hailing fleets. It enables faster recharge compared to even the fastest chargers, which can substantially reduce charging downtime, improve fleet utilization, and thereby enhance business sustainability. In this regard, battery swapping is a competitive solution for large commercial fleets and is attracting increasing attention from the industry recently. For instance, BP has partnered with Chinese battery swapping specialist Aulton to offer services for taxis, ride-hailing vehicles, and other passenger cars \cite{chinadaily}. San Francisco-based startup Ample provides battery swapping services for Uber and other commercial EVs using modular power packs. They estimate that a ride-hailing driver only needs to spend 3.6 hours per week charging their EV using battery swapping, which is 86\% less time than using fast charging \cite{ample}. Nevertheless, we must recognize that scaling a battery swapping network requires a longer timeline due to the higher deployment costs.

To capitalize on the advantages of plug-in charging and battery swapping, we propose a {\em multimodal charging network} in which charging stations and battery swapping stations are jointly built to support a ride-hailing EV fleet. Our central thesis is predicated on the observation that charging stations are cost-effective, making them ideally suitable for scaling up EV adoption in the beginning, while battery swapping stations offer quicker turnaround and can be deployed in tandem with charging stations to improve fleet utilization and reduce operational costs. To realize this vision, we consider a multi-stage infrastructure planning problem in which a profit-maximizing platform jointly deploys charging and battery swapping stations with a limited budget and simultaneously provides ride-hailing services over a transportation network. An optimization model is proposed to decide when and where to deploy each type of charging facilities, along with the ride-hailing operational decisions such as pricing strategies. The key contributions of this paper are summarized as follows:
\begin{enumerate}
\item This paper examines the benefits of a multimodal charging network consisting of both charging stations and battery swapping stations, with an emphasis on exploring their synergistic benefits. Despite extensive research on charging infrastructure planning, the complementary effects between these two facilities have been overlooked in the existing literature. To the best of our knowledge, this paper makes the first attempt to investigate the multimodal charging network for ride-hailing EV fleets by considering the joint deployment of charging stations and battery swapping stations.
    
\item We propose a multi-stage charging network expansion model to capture the coupled infrastructure planning and operational decisions in a ride-hailing network. This model integrates fundamental elements, including demand elasticity, drivers' charging behavior, charging and swapping waiting times, and their dependence on fleet status and charging infrastructure. Specifically, the stationary charging demand is analytically derived given passenger demand patterns and vehicle flows. The congestion at different facilities is captured by queuing-theoretic models, and drivers' decisions on charging location and mode are characterized by a spatial charging equilibrium model. 

\item We establish a tight upper bound for the joint deployment problem despite its nonconcavity. A decomposable problem is developed through proper model relaxations. By leveraging the favorable structures of the relaxed problem, we are able to obtain an accurate estimation of the globally optimal solution to the original problem, enabling us to verify the optimality of the solution obtained. We show that our approach provides a high-quality upper bound with an optimality gap of around 3\%.

\item Through extensive numerical studies, we evaluate the effectiveness of joint deployment and gain valuable insights into infrastructure planning. We find that joint planning of charging and battery swapping stations can harness the complementary benefits of the two facilities, leading to improved profit and reduced charging costs. In particular, based on a case study for Manhattan, New York City, the total profit of joint planning is 11.7\% higher than swapping-only deployment under a tight budget, and 17.5\% higher than charging-only deployment under a generous budget. Furthermore, compared to charging-only deployment, joint planning yields a charging cost reduction of 18.8\% under a limited budget and 44.4\% under a sufficient budget. However, a multimodal charging network may result in a higher volume of cross-zone traffic than a unimodal charging network due to driver's charging behaviors, which reveals the potential impacts of infrastructure planning on urban traffic congestion and management.
\end{enumerate}

\section{Literature Review}
\subsection{Modeling Electric Ride-Hailing Markets}
In ride-hailing markets, EV fleets differ from the gasoline-powered one in that EV charging would intermittently interrupt services and result in lengthy downtime, which impacts the platform's operation and necessitates special consideration in the formulation of mathematical models, the design of operational strategies, as well as the formation of regulatory policies.

One of the key topics for electric ride-hailing markets is the charging and operational strategies for EVs. Many studies focus on capturing the distinct driving patterns of ride-hailing EV drivers. For instance, Ke et al. \cite{ke_modelling_2019} modeled drivers' behavior in a ride-sourcing market with EVs and gasoline vehicles, developed a time-expanded network to sketch out drivers' working and recharging schedules, and found that passengers, drivers, and the platform are all worse off when the charging price increases. Qin et al. \cite{qin_development_2020} used a multi-state model to identify two distinct driving patterns for ride-hailing EVs, which outline when, for how long, and under what state of charge a driver will decide to recharge an EV. Alam \cite{alam_optimization-based_2022} emulated the daily trip patterns of ride-sourcing EVs using an optimization-based method and investigated the induced charging demand using agent-based simulation. Besides, a growing body of literature highlights the operations of electric ride-hailing markets. Shi et al. \cite{shi_operating_2020} formulated a dynamic vehicle routing problem (VRP) to capture the operations of a ride-hailing EV fleet and solved the problem using reinforcement learning. Kullman et al. \cite{kullman_dynamic_2022} also employed deep reinforcement learning to develop operational polices for ride-hailing EVs, including vehicle assignment, charging, and repositioning decisions. Pricing strategies are considered in \cite{hong_optimal_2022}, where Hong and Liu examined the optimal price and wage of green ride services while taking into account the opportunity cost of EV drivers as well as passengers' willingness to pay for eco-friendly trips. Ding et al. \cite{ding_optimal_2022} considered the integration of ride-hailing and vehicle-to-grid (V2G) services and developed a game-theoretic framework to characterize the market equilibrium. Cai et al. \cite{cai_integrating_2023} investigated the operations of integrated ride-sourcing and EV charging services in a nested two-sided market framework. Furthermore, related research is expanded into shared autonomous electric vehicles (SAEVs) \cite{chen_operations_2016}\cite{turan_dynamic_2020}\cite{yi_framework_2021}\cite{al-kanj_approximate_2020}. Among others, Turan et al. \cite{turan_dynamic_2020} considered pricing, routing, and charging strategies. Yi and Smart \cite{yi_framework_2021} focused on joint dispatching and charging management. Al-Kanj et al. \cite{al-kanj_approximate_2020} used approximate dynamic programming to develop high-quality pricing, dispatching, repositioning, charging, and fleet sizing strategies for the operations of an SAEV fleet.

Another important topic is on policies aimed at promoting the electrification of the ride-hailing industry. Different policy options have been explored in the literature, with a primary focus on pricing incentives that includes subsidy for EV purchases,  incentive schemes for EV rental, and financial support to infrastructure supply \cite{fleming_policy_2020}\cite{liu_regulatory_2022}\cite{slowik_how_2019}\cite{mo_modeling_2020}\cite{hall_guide_2021}\cite{lazer_electrifying_2021}. Particularly, Liu et al. \cite{liu_regulatory_2022} examined the market response and electrification levels under three different regulatory policies, including annual permit fees, differential trip-based fees, and differential commission caps. They indicated that the last policy is the most cost-efficient and can simultaneously benefit drivers and passengers. Slowik et al. \cite{slowik_how_2019} also found that well-designed taxes and fee structures can make EVs the most economically attractive technology for ride-hailing fleets. As per \cite{lazer_electrifying_2021}, governments should shift the target of current subsidies to the most intensively-used vehicles and the people most in need of financial support, for instance, ride-hailing drivers. And ensuring access to overnight charging is another crucial step to dismantling barriers to vehicle electrification. However, Mo et al. \cite{mo_modeling_2020} showed that given a limited budget, governments should subsidize charging infrastructure before supporting EV purchase. In summary, all aforementioned works focus on the operational strategies and/or public policies for electric ride-hailing fleets, while the planning of charging infrastructure for the ride-hailing network is not explicitly considered.

\subsection{Infrastructure Planning for EV fleets}
Charging infrastructure is a crucial component in a green and sustainable mobility ecosystem. In the existing literature, the planning strategies of EV infrastructure, including both charging stations and battery swapping stations, have been extensively examined (see \cite{shen_optimization_2019} for a comprehensive summary). Unlike private EVs that typically charge at home, commercial EV fleets have significantly higher mileage and rely heavily on public charging stations, requiring extra infrastructure planning considerations. In light of this, we review relevant studies that focus on infrastructure planning for commercial EV fleets, such as those in electric ride-hailing markets.

A recurrent topic in the literature is planning charging stations for EV fleets that offer mobility-on-demand services over the transportation network. Brandstätter et al. \cite{brandstatter_determining_2017} studied optimal locations of charging stations for a car-sharing EV fleet considering demand stochasticity. Roni et al. \cite{roni_optimal_2019} tackled charging station allocation and car-sharing EV assignment in a joint optimization framework and confirmed that charging time is a dominant component of the total fleet downtime. Ma and Xie \cite{ma_optimal_2021} also jointly considered infrastructure planning and vehicle assignment problem for ride-hailing EVs and formulated a bi-level problem to minimize the charging operation time. Alam and Guo \cite{alam_charging_2022} took drivers' value of time into account in planning charging stations for ride-hailing EVs. Bauer et al. \cite{bauer_electrifying_2019} adopted agent-based simulation to analyze charging infrastructure for electrified ride-hailing fleets and suggested that the cost of charging infrastructure accounts for a small fraction of the total cost of an electric ride-hailing trip. Following Lyft's announcement to reach 100\% EVs, Jenn \cite{jenn_charging_2021} examined the infrastructure deployment necessary to meet charging demand from fully-electrified fleets in California and found that a shift from daytime charging to overnight charging can significantly reduce costs.

The planning of battery swapping stations also attracts a lot of research focus. Mak et al. \cite{mak_infrastructure_2013} developed a robust optimization framework to deploy battery swapping facilities and investigated how battery standardization and technology advancements will impact the optimal deployment strategy. Targeted at refueling electric taxicabs, Wang et al. \cite{wang_toward_2018} jointly tackled the deployment of battery swapping stations and charging scheduling of EV fleets and designed a planning scheme to decrease in-station queuing time and improve quality of service of the taxicab fleet. Infante et al. \cite{infante_optimal_2020} considered stochastic charging demand and formulated a two-stage optimization with recourse, where battery swapping stations are deployed in the planning stage and the infrastructure network is coordinated in the operation stage. Further, Yang et al. \cite{yang_deploying_2021} proposed a data-driven approach to deploying and operating battery swapping stations with centralized management. Moon et al. \cite{moon_locating_2020} examined the locating problem of battery swapping stations in a passenger traffic network to achieve the seamless operation of electric buses.

Despite a large volume of literature on charging infrastructure planning for electric ride-hailing fleets, most existing works focus on either charging or battery swapping stations, without considering their complementary effects. Among the few early attempts, Zhang et al. \cite{zhang_joint_2021} studied the joint planning of swapping and charging stations for private EVs, and Zhang et al. \cite{zhang_towards_2020} explored the charging demand management for electric taxis in the presence of a hybrid infrastructure network. In contrast, this paper distinguishes itself from all previous studies in two key aspects. Firstly, we consider a multimodal charging network where charging stations and battery swapping stations are jointly deployed to overcome their respective limitations and elicit the synergy value. Secondly, we integrate infrastructure planning decisions and operational decisions of the ride-hailing platform within a unified framework to account for their interdependence. To the best of our knowledge, these considerations have not been studied in the literature.

\section{Mathematical Formulation}
We consider a platform that offers on-demand ride-hailing services using an EV fleet over a transportation network. To facilitate EV charging, the platform collaborates with charging infrastructure providers to devise a long-term infrastructure deployment strategy spanning $T$ stages,\footnote{This is to emulate a coalition of ride-hailing platforms and infrastructure providers, which follows the industry practice. For example,  Didi and BP also join forces to build an EV charging network in China \cite{BP_and_Didi}. Charging network provider Wallbox partners with Uber to install charging stations for Uber drivers \cite{wallbox}. Ample teams up with Sally, an EV fleet operator, to deploy swapping stations for ride-hailing and other service vehicles in metropolitan cities in the US \cite{ample_and_sally}.} each corresponding to a duration of 1-2 year, while subject to a limited budget.\footnote{The budget may come from government financial subsidies aimed at promoting charging infrastructure development, which is common in various jurisdictions \cite{protocal}\cite{EV_europe}.} At each stage, the platform determines how many charging stations and/or battery swapping stations to build, where to deploy charging infrastructure, and how to price ride-hailing services under the existing infrastructure network, etc. Note that infrastructure planning decisions are typically made on a timescale of years, whereas operational decisions, influenced by time-varying travel demand, impact market outcomes on a much shorter timescale, ranging from minutes to hours. In this light, the steady-state equilibrium at each stage is assumed to inform infrastructure planning. This assumption is commonly made in many existing works on ride-hailing markets  \cite{liu_regulatory_2022}\cite{zhou_competition_2022}\cite{ma_optimal_2021}. The remainder of this section details a mathematical model for the optimal deployment of a multimodal charging network.

\subsection{Charging Network Expansion}
The platform deploys a multimodal charging network in a city consisting of multiple zones. Let $x^k_{i,t}$ denote the number of charging/swapping stations built in zone $i$ up to stage $t$, where $k\in\m{K}=\{c,s\}$ is the index of facility mode, with $c$ and $s$ representing charging and swapping stations, respectively. Note that the planning decisions $x^k_{i,t}$ are considered continuous in this study, which suffices to evaluate the effectiveness of joint planning while also simplifying the solution process. Infrastructure deployment is subject to the budgetary constraint:
\begin{equation}\label{eqn:budget}
    \sum_{i\in\m{I}}\sum_{k\in\m{K}} \gamma_k x^k_{i,t} \leq b_t,
\end{equation}
where $\gamma_k$ is the unit deployment costs for different facilities, $b_t$ is the cumulative budget as of stage $t$, $\m{I}$ is the set of zones, and $\m{T}=\{1,\cdots,T\}$ is the set of stages within the planning horizon. This constraint implies that any unused budget during stage $t$ can be allocated to future stages. Define $x^k_{i,0}$ as the initial number of facilities. The expansion of charging network complies the following relation:
\begin{equation}\label{eqn:expansion}
    x^k_{i,t} - x^k_{i,t-1} \geq 0,
\end{equation}
indicating that the charging infrastructure built at each stage will be available for future use.

Given the infrastructure network at each stage $t$, the platform determines its operational strategies for profit maximization. These decisions involve both the supply and demand sides of ride-hailing markets, which are intimately coupled with infrastructure planning. To characterize this intrinsic correlation, we observe that the impacts of ride-hailing operations can be analyzed in the following sequence: (a) the platform set service prices over the transportation network to steer passenger demand; (b) passenger demand induces a flow pattern of EVs and generates a spatial distribution of charging demand; (c) drivers strategically decide charging locations and modes to minimize their individual charging downtime, leading to an equilibrium stage; (d) drivers' charging behavior, in turn, affects the total downtime of the entire fleet, thereby influencing infrastructure planning decisions. In subsequent subsections, we will present the platform's operational decisions following the aforementioned sequence.

\subsection{Passenger Demand Management}
Ride-hailing passengers are sensitive to price and waiting time. We denote the ride fare for a trip from zone $i$ to $j$ at stage $t$ as $p_{ij,t}$, and the average waiting time before each passenger is picked up as $w^p_{i,t}$. The arrival rate of passengers can be represented as:
\begin{equation}\label{eqn:demand model}
    q_{ij,t} = \bar{q}_{ij} F_p\left(p_{ij,t} + \beta w^p_{i,t}\right),
\end{equation}
where $\bar{q}_{ij}$ is the total travel demand from zone $i$ to $j$, including all passengers regardless of travel modes; $\beta$ is passengers' value of time; and $F_p(\cdot)$ is a strictly decreasing function representing the proportion of travelers using ride-hailing services. Since $F_p(\cdot)$ is strictly decreasing, it is invertible and the price $p_{ij,t}$ can be expressed as a function of $q_{ij,t}$ and $w^p_{i,t}$. Note that for a given $w^p_{i,t}$, the revenue, $p_{ij,t} q_{ij,t}$, is known to be concave in the demand rate $q_{ij,t}/\bar{q}_{ij}$ for various demand functions, such as linear, exponential, logit, and others \cite{chen_real-time_2023}. This property will be leveraged to facilitate the derivation of our proposed upper bounds. We therefore use the demand for each OD pair $q_{ij,t}$ as a decision variable in place of the price $p_{ij,t}$ hereinafter.

Assuming a first-come-first-serve matching principle, where the platform promptly assigns each rider to the nearest available vehicle, Passengers' waiting time only depends on the number of idle vehicles, denoted as $N_{i,t}^{v}$. This waiting time $w_{i,t}^p$ can be captured by the well-established square root law \cite{arnott_taxi_1996}\cite{li_regulating_2019}:
\begin{equation}\label{eqn:pickup time}
    w^p_{i,t} = \dfrac{\phi}{\sqrt{N^{v}_{i,t}}},
\end{equation}
where $\phi$ is a model parameter depending on the geometry of the city. This functional form is intuitive. As the number of idle vehicles increases subject to a fixed spatial distribution, the distance between a passenger and each of the vehicles will decrease in proportion to the square root of the number of idle vehicles. A rigorous proof for the case of uniform distribution can be found in \cite{li_regulating_2019}. Note that the first-come-first-serve matching mechanism is a common assumption in the existing literature on ride-hailing markets \cite{castillo_surge_2017}\cite{ke_pricing_2020}\cite{mo_modeling_2020}\cite{zhou_competition_2022}. We believe that it does not limit the model's generality since as the proposed methodology in this paper does not hinge on the square root structure of \eqref{eqn:pickup time}.

Due to the asymmetric passenger demand, the platform further dispatches vehicles from one zone to another to relocate the fleet. Let $r_{ij,t}$ denote the rebalancing flow from zone $i$ to $j$. In the long run, the vehicle inflow and outflow in each zone should be equal, leading to the following flow balance constraint:
\begin{equation}\label{eqn:flow balance}
    \sum_{j\in\m{I}} (q_{ij,t} + r_{ij,t}) = \sum_{j\in\m{I}} (q_{ji,t} + r_{ji,t}),
\end{equation}
where the left-hand side represents the total flow leaving zone $i$, and the right-hand side denotes the total flow entering zone $i$.

\subsection{Potential Charging Demand}
The stationary distribution of charging demand refers to the long-run probability of an EV requiring recharging in each zone. In this study, we assume that each driver will recharge the EV when its state of charge falls below a random threshold, with a mean value $\varepsilon$. As a result, a certain proportion of EVs will run out of power in each zone in the long run. Note that different EV drivers may adopt different charging schedules, characterized by distinct threshold values. However, when examining the stationary charging demand, we only need to focus on the average value of this threshold. We further assume that prior to reaching the threshold, the movement of an EV across the transportation network is completely governed by random passenger demand and is not affected by the specific threshold value. Consequently, the spatial distribution of charging demand aligns with the equilibrium fleet distribution, which, in turn, depends on the movement of the EV fleet and is inherently influenced by the platform's operational decisions.

Following the idea in \cite{qian_stationary_2019}, we consider the EV fleet movement as a random walk with a certain transition matrix $P$. To derive the transition matrix, we first note that the operating hours of EVs consist of idle trips (waiting for the next passenger), pickup trips (heading to pick up an assigned passenger), delivery trips (transporting passengers to destinations), and rebalancing trips (relocating from one zone to another with empty seats). The first two categories are regarded as intra-zone trips, while the last two types can be either inter-zone or intra-zone trips depending on destinations. In this light, we define $D_{ij,t}$ to capture the average number of EVs that travel from zone $i$ to zone $j$ (either with or without a passenger):
\begin{subnumcases}
    {D_{ij,t} = \label{eqn:demand matrix}}
    N^{v}_{i,t} + \sum_{j^\prime\in\m{I}} q_{ij^\prime,t} w^p_{i,t} + q_{ij,t} t_{ij} + r_{ij,t} t_{ij}, \quad & $i = j$, \label{eqn:demand matrix intra} \\
    q_{ij,t} t_{ij} + r_{ij,t} t_{ij}, \quad & $i \neq j$. \label{eqn:demand matrix inter}
\end{subnumcases}
The last three terms in \eqref{eqn:demand matrix intra} follow the Little's law and respectively represent the number of EVs undergoing pickup, delivery, and rebalancing trips. With this information, we can compute the transition matrix $P_t$:
\begin{equation}\label{eqn:transition matrix}
    P_{ij,t} = \dfrac{D_{ij,t}}{\sum_{j\in\m{I}} D_{ij,t}},
\end{equation}
which captures the probability of an EV departing for zone $j$ given that it is currently in zone $i$. Let $n_{i,t} > 0$ denote the probability of an EV staying in each zone at equilibrium. Based on the above analysis, the row vector $n_t = [n_{i,t}]_{i\in\m{I}}$ further represents the stationary charging demand distribution and thus satisfies $\sum_{i\in\m{I}} n_{i,t} = 1$. The derivation of this distribution follows the proposition below:
\begin{proposition}\label{prop:n existence}
    For any $q_{ij,t} > 0,\forall i,j\in\m{I}$ and the corresponding transition matrix $P_t = [P_{ij,t}]_{i,j\in\m{I}}$, the stationary charging demand distribution $n_t$ exists and is uniquely given by
    \begin{equation}\label{eqn:demand distribution}
        n_t P_t = n_t.
    \end{equation}
\end{proposition}
\begin{proof}
    Since $q_{ij,t} > 0$, we know that $P_t$ is element-wise strictly positive. Additionally, every row of $P_t$ adds to 1, so it is a row stochastic matrix. According to the Perron–Frobenius theorem \cite{strang_introduction_1993}, $P_t$ has an eigenvalue 1 and the corresponding left eigenvector is unique and has all-positive entries. As a result, solving \eqref{eqn:demand distribution} yields a unique stationary charging demand distribution $n_t$.
\end{proof}

The assumption of strict positivity for $q_{ij,t}$ is a mild one as it only requires the platform to provide services for each OD pair. Given $n_t$, we further define $\Lambda_{i,t}$ as the potential charging demand in each zone induced by EV movement, which clearly satisfies:
\begin{equation}\label{eqn:scaling n}
    \Lambda_{i,t} = n_{i,t} \sum_{j\in\m{I}} \Lambda_{j,t}.
\end{equation}
This potential charging demand $\Lambda_{i,t}$ can be accommodated by charging or battery swapping stations located in different zones. When an EV needs to recharge, the driver decides the charging location and mode. Let $f^k_{ij,t}$ denote the EV flow traveling from zone $i$ to $j$ to recharge using mode $k$ at stage $t$. Consequently, the potential charging demand $\Lambda_{i,t}$ can be expressed as the sum of these flows:
\begin{equation}\label{eqn:potential charging demand}
    \Lambda_{i,t} = \sum_{j\in\m{I}}\sum_{k\in\m{K}} f^k_{ij,t}.
\end{equation}
Let $\omega$ denote the average rate of electricity consumption, and let $\Omega$ denote the average battery capacity. At equilibrium, the electricity consumed by the entire EV fleet should match the total energy replenished through plug-in charging and battery swapping. This leads to the following energy balance condition: 
\begin{equation}\label{eqn:energy balance}
    \omega \left[ \sum_{i\in\m{I}} N^v_{i,t} + \sum_{i\in\m{I}}\sum_{j\in\m{I}} q_{ij,t} w^p_{i,t} + \sum_{i\in\m{I}}\sum_{j\in\m{I}} q_{ij,t} t_{ij} + \sum_{i\in\m{I}}\sum_{j\in\m{I}} r_{ij,t} t_{ij} + \sum_{i\in\m{I}}\sum_{j\in\m{I}}\sum_{k\in\m{K}} f^k_{ij,t} t_{ij} \right] = (1-\varepsilon)\Omega \sum_{i\in\m{I}} \Lambda_{i,t}
\end{equation}
where $t_{ij}$ is the travel time from zone $i$ to $j$ and $\varepsilon$ is the average energy threshold below which an EV will top up the battery. In \eqref{eqn:energy balance}, the left-hand side is derived from Little's law and represents the total energy consumption, while the right-hand side represents the power input to the entire fleet. The terms on the left-hand side account for the energy consumed by EVs that are idling, picking up passengers, delivering passengers, relocating, and en-route to recharge, respectively.

So far, we have characterized how the platform's operational decisions affect the potential charging demand. Specifically, \eqref{eqn:demand distribution} reveals the equilibrium spatial distribution of charging demand, while \eqref{eqn:energy balance} provides a scaling factor for \eqref{eqn:demand distribution} to ensure a balanced energy input and output. It is important to note that the potential charging demand is endogenous in the sense that it depends on the transition matrix $P_t$, which further depends on the platform's operational decisions.

\begin{remark}
    {Ideally, charging flows $f^k_{ij,t}$ are part of the vehicle movement and can impact the equilibrium fleet distribution. Therefore, they should be included in \eqref{eqn:demand matrix} when deriving the transition matrix $P_t$. However, incorporating charging flows in this way would result in a more complex model with compromised tractability. As a result, we have decided to simplify the model by assuming that the transition matrix is independent of charging flows. We acknowledge that this simplification is a limitation of our proposed model and recognize the potential for some deviation in the results. However, we argue that this deviation turns out to be insignificant since charging flows typically represent a small fraction of the total EV movement. In a realistic simulation study, we find that the deviation of charging demand distributions due to the neglection of charging flows $f^k_{ij,t}$ in \eqref{eqn:demand matrix} is less than 0.1\% (measured by L2 distance). Hence, we leave the development of a more comprehensive model that incorporates charging flows in the transition matrix as an avenue for future research.}
\end{remark}

\subsection{Spatial Charging Equilibrium}
Given the potential charging demand, drivers make strategic decisions regarding where to recharge and which charging mode to use so as to minimize their charging costs. The aggregate choices of drivers, characterized by the charging flow $f^k_{ij,t}$, collectively establish a spatial charging equilibrium across different zones. To formally capture this equilibrium, we first introduce the notion of actual charging demand:
\begin{equation}\label{eqn:actual charging demand}
    \lambda^k_{i,t} = \sum_{j\in\m{I}} f^k_{ji,t},
\end{equation}
which defines the total demand for charging facilities of mode $k$ in zone $i$ at stage $t$. In addition, the charging cost of each EV driver consists of the travel time to the selected zone, the waiting time at the visited facility, and the service time of the chosen charging mode:
\begin{equation}
    c^k_{ij,t} = t_{ij} + t_k + w^k_{j,t}\left(x_{j,t},\lambda^k_{j,t}\right).
\end{equation}
Here, $t_k$ corresponds to the service time associated with mode $k$, and $w^k_{j,t}$ denotes the queueing time, which is a function of facility supply $x_{j,t}$ and charging demand $\lambda^k_{j,t}$. The precise relationship between $w^k_{j,t}$ and $x_{j,t}$, as well as $\lambda^k_{j,t}$, will be elaborated upon in the subsequent subsection. We then introduce the definition of spatial charging equilibrium:
\begin{definition}[Spatial Charging Equilibrium]
    The system is said to reach spatial charging equilibrium when no EV driver from the same zone can reduce her charging cost by unilaterally deviating from the selected charging location or mode.
\end{definition}
Based on this definition, the equilibrium conditions can be equivalently expressed as
\begin{subnumcases}{\label{eqn:equilibrium condition}}
    f^k_{ij,t} \geq 0, \\
    c^k_{ij,t} - u_{i,t} \geq 0, \\
    f^k_{ij,t} (c^k_{ij,t} - u_{i,t}) = 0, \label{eqn:complementarity}
\end{subnumcases}
together with the flow conservation constraint \eqref{eqn:potential charging demand}, where $u_{i,t}$ denotes the equilibrium charging cost for each EV from zone $i$ at stage $t$. These conditions state that if an EV travels from zone $i$ to $j$ to recharge using mode $k$, the associated charging cost is equal to the equilibrium cost for zone $i$, while those not selected location and mode will incur a charging cost greater than or at least equal to the equilibrium cost $u_{i,t}$. Thereby, if the flow pattern $f^k_{ij,t}$ satisfies these conditions, no drivers can be better off by unilaterally changing their decisions.

\subsection{Queueing Time at Stations}
This subsection will detail the waiting times at charging and battery swapping stations. We proceed by first assuming that different facilities are uniformly distributed within each zone, and that EV drivers only strategically determine the zone and mode for recharging, but will always visit the closest facility of the selected mode. Therefore, the average arrival rate of EVs for each facility of mode $k$ is $\bar{\lambda}^k_{i,t} = \lambda^k_{i,t}/x^k_{i,t}$.

{\em Queueing time at charging stations:} Assuming the arrival of EVs follows a Poisson distribution with rate $\bar{\lambda}^c_{i,t}$ and the charging time follows an exponential distribution with an average value of $t_c$, we can model the charging process as an $M/M/C/K$ queue. Here, $C$ represents the capacity of each charging station (i.e., the number of chargers), and $K$ represents the queue capacity (i.e., only a maximum of $K$ customers are allowed to stay in the queue). The mean queue length for this queueing system is well-established \cite{hillier_introduction_2004} and can be expressed as:
\begin{equation}
    L_c = \dfrac{\rho P_0 (\bar{\lambda}^c_{i,t} t_c)^C }{C!(1-\rho)^2} \left[ 1 - \rho^{K-C+1} - (1-\rho)(K-C+1)\rho^{K-C} \right],
\end{equation}
where $\rho = \bar{\lambda}^c_{i,t} t_c/C$, and
\begin{subnumcases}
    {P_0 = }
    \left[ \dfrac{(\bar{\lambda}^c_{i,t} t_c)^C}{C!} \dfrac{1-\rho^{K-C+1}}{1-\rho} + \sum\limits_{a=0}^{C-1} \dfrac{(\bar{\lambda}^c_{i,t} t_c)^a}{a!}\right]^{-1}, & $\t{if } \rho \neq 1$, \\
    \left[ \dfrac{(\bar{\lambda}^c_{i,t} t_c)^C}{C!} (K-C+1) + \sum\limits_{a=0}^{C-1} \dfrac{(\bar{\lambda}^c_{i,t} t_c)^a}{a!}\right]^{-1}, & $\t{if } \rho = 1$.
\end{subnumcases}
Note that if $\rho = 1$, L'Hôpital's rule should be applied twice to the expression of $L_c$. Due to the limited queue capacity, new arrivals will be blocked if the queue is full. The blocking rate is given by
\begin{equation}
    P_K = \dfrac{(\bar{\lambda}^c_{i,t} t_c)^K}{C! C^{K-C}} P_0.
\end{equation}
The average waiting time at each charging station can thus be obtained by applying Little's law:
\begin{equation}\label{eqn:charging queue}
    w^c_{i,t} = \dfrac{L_c}{\bar{\lambda}^c_{i,t} (1 - P_K)}.
\end{equation}

\begin{figure}
    \centering
    \includegraphics[width=0.8\linewidth]{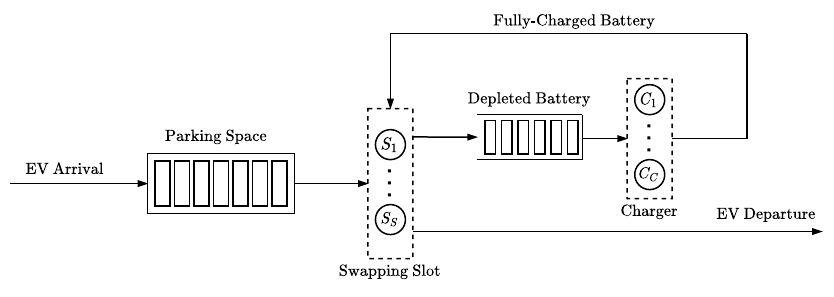}
    \caption{The mixed queueing network for a battery swapping station.}
    \label{fig:queueing network}
\end{figure}

{\em Queueing time at swapping stations:} Compared to charging stations,
battery swapping stations have much more complex dynamics. If fully-charged batteries are available when an EV arrives, it will be immediately maneuvered into a changeover bay for battery replacement. Otherwise, EVs will queue up for the next fully-charged power pack as long as the queue is not full upon arrival. The depleted batteries will be plugged in right after being removed from EVs. To model these dynamics, we adopt the framework in \cite{tan_queueing_2014} and formulate the system as a mixed queueing network. As illustrated in \autoref{fig:queueing network}, the queueing network includes an open queue for EVs and a closed queue for batteries. The former captures the arrival and departure of EVs, while the latter models the circulation of batteries. The battery queue is considered closed because the arrival of depleted batteries and the departure of fully-charged batteries are in one-to-one correspondence. Consequently, each station maintains a fixed number of batteries in circulation.

To ensure consistency and fair comparison between the two facilities, we consider each swapping station is equipped with $C$ chargers for recharging depleted batteries and has a queue capacity of $K$. Additionally, we denote $S$ as the number of automated swapping slots, $B$ as the number of batteries in circulation, and $t_s$ as the constant service time for each swap. The arrival of EVs at each swapping station is also assumed to be Poisson with rate $\bar{\lambda}^s_{i,t}$. Define $g(u)$ as the probability of having $u$ new-coming EVs during any time interval $(t, t+t_s]$. Based on the probability mass function of Poisson distribution, we derive $g(u)$ as follows:
\begin{equation}\label{eqn:arrival}
    g(u) = \dfrac{(\bar{\lambda}^s_{i,t} t_s)^u \exp(-\bar{\lambda}^s_{i,t} t_s)}{u!}, \, u \geq 0.
\end{equation}
Further let $f(u,v)$ represent the probability of having $v$ batteries complete charging during $(t, t+t_s]$ provided that there are $u$ batteries being charged at time $t$. Clearly, $f(u,v)$ is supported on $u \in [0,C]$ and $v \in [0,u]$. As aforementioned in the charging queue, suppose the battery charging time follows an independent and identical exponential distribution with mean $t_c$. Then $f(u,v)$ is equivalent to having $v$ ``success'' in $u$ Bernoulli trials, each of which has a success probability $1 - \exp(-t_s/t_c)$:
\begin{equation}\label{eqn:battery circulation}
    f(u,v) = \begin{pmatrix} u \\ v \end{pmatrix} \left(1 - \exp(-t_s/t_c)\right)^v \exp(-(u-v)t_s/t_c).
\end{equation}
With \eqref{eqn:arrival} and \eqref{eqn:battery circulation} respectively capturing the dynamics of EV arrival and battery circulation, we can proceed to establish the queueing equilibrium. The system can be characterized by a state tuple $(i,j)$ representing that there are $i$ EVs (including those waiting and being served) and $j$ fully-charged batteries in the station. Denote $G_t(i,j)$ as the probability of system state $(i,j)$ at time $t$ and denote $G_{t+t_s}(\hat{i},\hat{j})$ as the probability of system state $(\hat{i},\hat{j})$ at time $t+t_s$. We can then compute $G_{t+t_s}(\hat{i},\hat{j})$ given $G_t(i,j)$:
\begin{subnumcases}
    {G_{t+t_s}(\hat{i},\hat{j}) = \label{eqn:steady state}}
    \sum_{i=0}^K \sum_{j=0}^B G_t(i,j) g(\hat{i}-i+\Delta) f(B-j,\hat{j}-j+\Delta), & $\hat{i} < K$, \\
    \sum_{i=0}^K \sum_{j=0}^B G_t(i,j) \left(1 - \sum_{d=1}^{K-1} g(d-i+\Delta) \right) f(B-j,\hat{j}-j+\Delta), & $\hat{i} = K$,
\end{subnumcases}
where $\Delta=\min\{i,j,S\}$ is the number of EVs leaving the station with fully-charged batteries during $(t,t+t_s]$. When the queue is not full at $t+t_s$, i.e., $\hat{i} < K$, the transition from $(i,j)$ to $(\hat{i},\hat{j})$ is equivalent to having $\hat{i}-i+\Delta$ EVs arriving and $\hat{j}-j+\Delta$ batteries finishing charging. However, if the queue capacity is reached, i.e., $\hat{i} = K$, there are at least $K-i-\Delta$ new arrivals during $(t,t+t_s]$, which happens with a probability of $1 - \sum_{d=1}^{K-1} g(d-i+\Delta)$. The equilibrium state distribution, denoted as $G(i,j) = \lim_{t \rightarrow \infty} G_t(i,j)$, is given by solving \eqref{eqn:steady state} together with the normalization equation $\sum_{i=0}^K \sum_{j=0}^B G(i,j) = 1$. Note that the equation system consists of $(K+1)\times(B+1)+1$ linear equations but only $(K+1)\times(B+1)$ variables, which may not necessarily have a unique solution. Whereas it is proven in \cite[Theorem 1]{tan_queueing_2014} that the corresponding embedded Markov chain is ergodic. Thereby we can always find a unique equilibrium state distribution $G(i,j)$, based on which we can derive the average waiting time $w^s_{i,t}$. An efficient approach to computing $G(i,j)$ without solving an equation system is detailed in \ref{appendix:swapping queue}. The expected queue length is $L = \sum_{i=0}^W \sum_{j=0}^B i \cdot G(i,j)$. Applying the Little's law yields
\begin{equation}\label{eqn:queue length}
    L_s = \bar{\lambda}^s_{i,t} (w^s_{i,t} + t_s) \left(1 - \sum_{j=0}^B G(K,j) \right),
\end{equation}
where $1 - \sum_{j=0}^B G(K,j)$ represents the probability that the queue is not full. Hence, the average waiting time at battery swapping stations is given by
\begin{equation}\label{eqn:swapping queue}
    w^s_{i,t} = \dfrac{\sum_{i=0}^K \sum_{j=0}^B i \cdot G(i,j)}{\bar{\lambda}^s_{i,t} \left(1 - \sum_{j=0}^B G(K,j) \right)} - t_s.
\end{equation}

In this study, charging and battery swapping queues are both considered finite. As a result, an EV would have to seek an alternative location and mode for recharging if the queue is full upon its arrival. However, it is important to note that we assume a large queue capacity in this study to prevent the case of a full queue. This consideration is based on the observation that drivers do not necessarily wait exactly at the stations, as these stations might have limited space to accommodate a long queue. Instead, drivers can wait in the vicinity of the charging facilities, such as roadside parking spaces or nearby parking lots, where a significantly larger number of vehicles can be held in queue. With a large queue capacity, the blocking rate is negligible for any EV arrival rate in a practical range. We show in \ref{appendix:queueing model} that the blocking rate becomes nonnegligible only when demand at each station is very high. In this case, the waiting time is prohibitively long, which deters EVs from using the facility at equilibrium. To avoid this, we further enforce the following constraint:
\begin{equation}\label{eqn:demand upper bound}
    0 \leq \lambda^k_{i,t} \leq \rho_k x^k_{i,t},
\end{equation}
which imposes an upper bound $\rho_k$ on the average arrival rate of each facility. This ensures the actual charging demand falls in a reasonable regime where the consequent waiting time remains tolerable (e.g., less than 1 hour) and the blocking rate remains marginal. From this perspective, constraint \eqref{eqn:demand upper bound} can be interpreted as a quality of service (QoS) guarantee for the infrastructure network.\footnote{Results in our case studies demonstrate that this constraint is always inactive, indicating that EV drivers will strategically avoid congested facilities with high chance of blockage.}
\begin{remark}
In situations where the queue capacity is small, the probability of blocking may become non-negligible, and the model has to be revised accordingly. One potential approach to address this concern is to consider a slightly different queuing model, where there are multiple charging stations from the same zone in each queue, and the driver would choose which station to approach within the same zone based on which station has the shortest queuing time. In this case, the blocking probability could be negligible, as the likelihood of all charging stations in the same zone being fully occupied is quite low. In this paper, we focus on the large capacity case and leave the other case as future work.
\end{remark}

\subsection{Decision-Making Problem}
The decision-making process in this study involves both the platform and EV drivers. The platform makes planning decisions ($x^k_{i,t}$) and operational decisions ($q_{ij,t},r_{ij,t},N^v_{i,t}$), based on which EV drivers decide their charging locations and modes through $f^k_{ij,t}$. These decisions not only influence endogenous variables such as $\Lambda_{i,t}$, $\lambda^k_{i,t}$, and $w^k_{i,t}$ but also determine the fleet size:
\begin{equation}\label{eqn:fleet size}
    N_t = \sum_{i\in\m{I}} N^v_{i,t} + \sum_{i\in\m{I}}\sum_{j\in\m{I}} q_{ij,t} w^p_{i,t} + \sum_{i\in\m{I}}\sum_{j\in\m{I}} q_{ij,t} t_{ij} + \sum_{i\in\m{I}}\sum_{j\in\m{I}} r_{ij,t} t_{ij}  + \sum_{i\in\m{I}}\sum_{j\in\m{I}}\sum_{k\in\m{K}} f^k_{ij,t} t_{ij} + \sum_{i\in\m{I}}\sum_{k\in\m{K}} \lambda^k_{i,t} (w^k_{i,t} + t_k)
\end{equation}
where the first four terms represent EVs in operations, while the last two terms account for the total number of EVs traveling to, waiting at, and being served in charging and battery swapping stations. The platform's revenue at each stage from ride-hailing services can be then represented as
\begin{equation}
    \Pi^r_t = \sum_{i\in\m{I}}\sum_{j\in\m{I}} p_{ij,t} q_{ij,t} - \gamma_e N_t - \gamma_p \sum_{i\in\m{I}}\sum_{j\in\m{I}}\sum_{k\in\m{K}} f^k_{ij,t} u_{i,t},
\end{equation}
where $\gamma_e$ is the unit cost of each EV and $\gamma_p$ is a penalty parameter. The last term penalizes the total time cost resulting from drivers' charging behaviors, which ensures that charging infrastructure is properly deployed to maintain charging costs within an acceptable range. On the other hand, the deployment cost of charging facilities at each stage is
\begin{equation}
    \Pi^c_t = \sum_{i\in\m{I}}\sum_{k\in\m{K}} \gamma_k (x^k_{i,t} - x^k_{i,t-1}).
\end{equation}
The platform stops building charging infrastructure after stage $T$. However, the infrastructure will stay and remain functional beyond time $T$ until it exceeds the lifespan, i.e., $H > T$. Therefore, the deployment cost is incurred in the first $T$ stages, while the platform's revenue and operational costs remain unchanged from stage $T+1$ onward until $H$. To reflect the preference for immediate costs or benefits over future ones and accounts for the time value of money and/or risk, we introduce a discount factor $\gamma \in (0,1)$ \cite{lin_multistage_2019}. The total profit over the entire operation horizon is then defined as
\begin{equation}
    \Pi = \sum_{t\in\m{T}} \gamma^{t-1} (\Pi^r_t - \Pi^c_t) + \dfrac{\gamma^T - \gamma^H}{1 - \gamma} \Pi^r_T.
\end{equation}
By equivalently replacing the driver's decision problem with equilibrium conditions \eqref{eqn:equilibrium condition}, we formulate the decision-making problem as a mathematical program with complementarity constraints (MPCC):
\begin{subequations}\label{eqn:original problem}
\begin{align}
    \underset{x,q,r,f,u,N^v}{\t{maximize}} \quad & \sum_{t\in\m{T}} \gamma^{t-1} (\Pi^r_t - \Pi^c_t) + \dfrac{\gamma^T - \gamma^H}{1 - \gamma} \Pi^r_T \label{eqn:original objective} \\
    \t{subject to} \quad & \t{constraints }\eqref{eqn:budget}-\eqref{eqn:fleet size} \notag \\
    & N^v_{i,t} > \hat{N}^v, \ \forall t\in\m{T},i\in\m{I}, \\
    & q_{ij,t} > 0, \ \forall t\in\m{T},i\in\m{I},j\in\m{I} \\
    & x^k_{i,0} \leq x^k_{i,t} \leq \hat{x}^k_{i}, \ \forall t\in\m{T},i\in\m{I},k\in\m{K},
\end{align}    
\end{subequations}
where $\hat{N}^v$ is the lower bound for the number of idle vehicle in each zone, which reflects the practice that the platform will maintain the vehicle supply above a certain threshold to ensure service quality; and $\hat{x}^k_{i}$ is the maximum number of facilities allowed to deploy in each zone because of practical restrictions, e.g., limited land supply. This optimization problem is inherently challenging to solve due to the complementarity constraints \eqref{eqn:equilibrium condition} and the nonconvex constraints \eqref{eqn:demand distribution}\eqref{eqn:charging queue}\eqref{eqn:swapping queue}.

\section{Solution Method}
In this section, we will present our approach for pinpointing stationary points for \eqref{eqn:original problem} and subsequently establish a theoretical upper bound based on model relaxation and decomposition.

\subsection{Finding Stationary Points}\label{sec:stationary point}
Nonlinear programming methods typically require constraint qualification to ensure the existence and boundedness of Lagrange multipliers and the satisfaction of Karush-Kuhn-Tucker (KKT) conditions. The complementarity constraints \eqref{eqn:equilibrium condition}, however, violate Mangasarian-Fromovitz Constraint Qualification (MFCQ) at any feasible point because of their combinatorial nature. Therefore, directly solving \eqref{eqn:original problem} using common nonlinear programming algorithms may not be able to obtain a local stationary point.

To address this issue, an auxiliary parameter $\epsilon$ is introduced to relax the complementarity condition \eqref{eqn:complementarity}:
\begin{equation}\label{eqn:relaxed complementarity}
    f^k_{ij,t} (c^k_{ij,t} - u_{i,t}) \leq \epsilon,
\end{equation}
which is a common relaxation scheme to iteratively compute the stationary point of MPCC \cite{scholtes_convergence_2001}\cite{ban_general_2006}\cite{chen_optimal_2020}. For any $\epsilon > 0$, the gradients of active constraints are linear independent. Therefore, the resulting optimization problem satisfies constraint qualification and can be solved using nonlinear programming algorithms. We can start from a relatively large value of $\epsilon$ and gradually decrease its value. At each iteration, we solve the corresponding optimization problem using the solution obtained from the previous iteration as an initial guess. The iterative process terminates at a sufficiently small $\epsilon$, where we obtain a solution that can be used to derive an exact solution to the original problem \eqref{eqn:original problem}.

\subsection{Derivation of Upper Bounds}
Nonlinear programming algorithms, although capable of finding local solutions, do not provide any performance guarantee. To address this limitation, our objective is to establish an upper bound for the original problem, which will enable us to verify the optimality of any obtained solution, even in the absence of concavity. However, it is important to note that finding a tight upper bound for nonconcave optimization problems is a challenging task due to the complexity and diversity of problem structures. Common approaches to establishing upper bounds heavily rely on relaxation techniques, such as convex relaxation and linear programming relaxation, which results in a simplified optimization problem that is easier to solve. In the same vein, we leverage the problem's characteristics and propose a relaxed formulation that presents favorable structures while preserving the essence of the original problem, based on which a high-quality upper bound can be found.

\subsubsection{Model Relaxation}
The original formulation is relaxed by assuming that the platform has direct control over the actual charging demand $\lambda^k_{i,t}$ in each zone by regulating vehicle flow $f^k_{ij,t}$. However, charging mode choice is still subject to equilibrium conditions, meaning that EV drivers only decide their charging modes within the zone assigned by the platform. This relaxation is grounded in the intuition that the actual charging demand is naturally close to the desired pattern of the platform, since the platform can induce the desired charging behavior through planning and operational decisions. Furthermore, the preferred charging patterns of drivers and the platform are partially aligned because of their mutual objective of reducing the charging cost. Given these considerations, we believe this relaxation will not result in a significant deviation.

After the relaxation, we can replace $f^k_{ij,t}$ with $\lambda^k_{i,t}$ in the formulation. Define the relaxed charging cost as:
\begin{equation}
    \tilde{c}^k_{i,t} = \tilde{t}_{i} + w^k_{i,t} + t_k,
\end{equation}
in which $\tilde{t}_i = \min_{\forall j}\{t_{ji}\}$ denotes the minimum time required to reach zone $i$. The equilibrium condition for charging mode choice then becomes:
\begin{subnumcases}{\label{eqn:relaxed equilibrium condition}}
    \lambda^k_{i,t} \geq 0, \\
    \tilde{c}^k_{i,t} - \tilde{u}_{i,t} \geq 0, \\
    \lambda^k_{i,t} (\tilde{c}^k_{i,t} - \tilde{u}_{i,t}) = 0,
\end{subnumcases}
where $\tilde{u}_{i,t}$ is the resulting equilibrium charging cost. By definition, we further have $\tilde{c}^k_{i,t} \leq c^k_{ji,t}$ and $\tilde{u}_{i,t} \leq u_{i,t}$. Consequently, the total charging cost can be expressed as:
\begin{equation}
    \sum_{i\in\m{I}}\sum_{k\in\m{K}} \lambda^k_{i,t} \tilde{u}_{i,t} = \sum_{i\in\m{I}}\sum_{j\in\m{I}}\sum_{k\in\m{K}} f^k_{ij,t} \tilde{c}^k_{j,t} \leq \sum_{i\in\m{I}}\sum_{j\in\m{I}}\sum_{k\in\m{K}} f^k_{ij,t} c^k_{ij,t} =
    \sum_{i\in\m{I}}\sum_{j\in\m{I}}\sum_{k\in\m{K}} f^k_{ij,t} u_{i,t},
\end{equation}
which is upper bounded by the original one. By replacing $t_{ij}$ in the second last term of \eqref{eqn:fleet size} with $\tilde{t}_i$, we can also relax the fleet size to be
\begin{equation}
\begin{aligned}
    \widetilde{N}_{t} &= \sum_{i\in\m{I}} N^v_{i,t} + \sum_{i\in\m{I}}\sum_{j\in\m{I}} q_{ij,t} w^p_{i,t} + \sum_{i\in\m{I}}\sum_{j\in\m{I}} q_{ij,t} t_{ij} + \sum_{i\in\m{I}}\sum_{j\in\m{I}} r_{ij,t} t_{ij} + \sum_{i\in\m{I}}\sum_{j\in\m{I}}\sum_{k\in\m{K}} f^k_{ij,t} \tilde{t}_j + \sum_{i\in\m{I}}\sum_{k\in\m{K}} \lambda^k_{i,t} (w^k_{i,t} + t_k) \\
    &= \sum_{i\in\m{I}} N^v_{i,t} + \sum_{i\in\m{I}}\sum_{j\in\m{I}} q_{ij,t} w^p_{i,t} + \sum_{i\in\m{I}}\sum_{j\in\m{I}} q_{ij,t} t_{ij} + \sum_{i\in\m{I}}\sum_{j\in\m{I}} r_{ij,t} t_{ij} + \sum_{i\in\m{I}}\sum_{k\in\m{K}} \lambda^k_{i,t} (\tilde{t}_i + w^k_{i,t} + t_k).
\end{aligned}
\end{equation}
Correspondingly, the energy balance condition \eqref{eqn:energy balance} becomes
\begin{equation}
    \sum_{i\in\m{I}} N^v_{i,t} + \sum_{i\in\m{I}}\sum_{j\in\m{I}} q_{ij,t} w^p_{i,t} + \sum_{i\in\m{I}}\sum_{j\in\m{I}} q_{ij,t} t_{ij} + \sum_{i\in\m{I}}\sum_{j\in\m{I}} r_{ij,t} t_{ij} = \sum_{i\in\m{I}}\sum_{k\in\m{K}} \lambda^k_{i,t} \left[\dfrac{(1-\varepsilon)\Omega}{\omega} - \tilde{t}_i\right].
\end{equation}
For notation brevity, we further rewrite the coefficients in \eqref{eqn:original objective} as $\sigma^c_t = \gamma^{t-1}$ and
\begin{subnumcases}
    {\sigma^r_t =}
    \gamma^{t-1}, & $t < T$, \\
    \dfrac{\gamma^{T-1} - \gamma^H}{1 - \gamma}, & $t = T$.
\end{subnumcases}
After these modifications, the original problem \eqref{eqn:original problem} is reformulated as follows:
\begin{subequations}\label{eqn:relaxed problem}
\begin{align}
    \underset{x,q,r,\lambda,\tilde{u},N^v}{\t{maximize}} \quad & \sum_{t\in\m{T}} \left(\sigma^r_t \Pi^r_t - \sigma^c_t \Pi^c_t \right) \\
    \t{subject to} \quad & \sum_{i\in\m{I}}\sum_{k\in\m{K}} \gamma_k x^k_{i,t} \leq b_t, \ \forall t\in\m{T}, \label{eqn:relaxed problem budget} \\
    & N^v_{i,t} > \hat{N}^v, \ \forall t\in\m{T},i\in\m{I}, \\
    & q_{ij,t} > 0, \ \forall t\in\m{T},i\in\m{I},j\in\m{I}, \\
    & x^k_{i,0} \leq x^k_{i,t} \leq \hat{x}^k_{i}, \ \forall t\in\m{T},i\in\m{I},k\in\m{K}, \\
    & x^k_{i,t} - x^k_{i,t-1} \geq 0, \ \forall t\in\m{T},i\in\m{I},k\in\m{K}, \label{eqn:relaxed expansion}\\
    & 0 \leq \lambda^k_{i,t} \leq \rho_k x^k_{i,t}, \ \forall t\in\m{T},i\in\m{I},k\in\m{K}, \\
    & \lambda^k_{i,t} (\tilde{c}^k_{i,t} - \tilde{u}_{i,t}) \leq \epsilon, \ \forall t\in\m{T},i\in\m{I},k\in\m{K}, \label{eqn:relaxed problem complementarity} \\
    & \tilde{c}^k_{i,t} - \tilde{u}_{i,t} \geq 0, \lambda^k_{i,t} \geq 0, \ \forall t\in\m{T},i\in\m{I},k\in\m{K}, \\
    & \sum_{j\in\m{I}} (q_{ij,t} + r_{ij,t}) = \sum_{j\in\m{I}} (q_{ji,t} + r_{ji,t}), \ \forall t\in\m{T},i\in\m{I}, \label{eqn:relaxed problem flow balance} \\
    & \sum_{i\in\m{I}} N^v_{i,t} + \sum_{i\in\m{I}}\sum_{j\in\m{I}} \left(q_{ij,t} w^p_{i,t} + q_{ij,t} t_{ij} + r_{ij,t} t_{ij}\right) = \sum_{i\in\m{I}}\sum_{k\in\m{K}} \lambda^k_{i,t} \left[\dfrac{(1-\varepsilon)\Omega}{\omega} - \tilde{t}_i\right], \ \forall t\in\m{T}. \label{eqn:relaxed problem energy balance}
\end{align}
\end{subequations}
Compared to \eqref{eqn:original problem}, the relaxed problem \eqref{eqn:relaxed problem} remains nonconcave but is more tractable owing to its decomposable structure, which can be leveraged to derive an upper bound subsequently.

\subsubsection{Model Decomposition}
The derivation of upper bounds is based on Lagrangian relaxation, whose performance is sensitive to the quality of Lagrange multipliers. Therefore, constraint \eqref{eqn:relaxed expansion} is removed from \eqref{eqn:relaxed problem} and the complementarity constraint \eqref{eqn:relaxed problem complementarity} is relaxed by $\epsilon$ as in \eqref{eqn:relaxed complementarity}. These modifications effectively improve the tractability of \eqref{eqn:relaxed problem} without compromising the validity of model relaxation. Let $\Phi=(\eta,\mu,\theta)$ be the Lagrange multipliers associated with \eqref{eqn:relaxed problem budget},\eqref{eqn:relaxed problem flow balance}, and \eqref{eqn:relaxed problem energy balance}, respectively. The partial Lagrangian is:
\begin{equation}
\begin{aligned}
    \m{L}(x,q,r,\lambda,\tilde{u},N^v)\Big|_{\Phi} =& \sum_{t\in\m{T}} \left(\sigma^r_t \Pi^r_t - \sigma^c_t \Pi^c_t \right) + \sum_{t\in\m{T}}\sum_{i\in\m{I}} \sum_{j\in\m{I}} (\mu_{i,t} - \mu_{j,t}) (q_{ij,t} + r_{ij,t}) \\
    & - \sum_{t\in\m{T}} \sum_{i\in\m{I}}\sum_{k\in\m{K}} \eta_t\gamma_k x^k_{i,t} + \sum_{t\in\m{T}}\sum_{i\in\m{I}} \theta_t N^v_{i,t} + \sum_{t\in\m{T}}\sum_{i\in\m{I}}\sum_{j\in\m{I}} \theta_t q_{ij,t} w^p_{i,t} \\
    & + \sum_{t\in\m{T}}\sum_{i\in\m{I}}\sum_{j\in\m{I}} \theta_t \left(q_{ij,t} + r_{ij,t}\right) t_{ij} - \sum_{t\in\m{T}}\sum_{i\in\m{I}}\sum_{k\in\m{K}} \theta_t \lambda^k_{i,t} \left[\dfrac{(1-\varepsilon)\Omega}{\omega} - \tilde{t}_i\right] 
\end{aligned}
\end{equation}
Since $w^k_{i,t}$ only depends on $x^k_{i,t}$ and $\lambda^k_{i,t}$, given the multipliers, the Lagrangian is separable in both spatial and temporal dimensions, i.e., $\m{L}(x,q,r,\lambda,\tilde{u},N^v) = \sum_{t\in\m{T}}\sum_{i\in\m{I}} \m{L}_{i,t}(x^k_{i,t},q_{ij,t},r_{ij,t},\lambda^k_{i,t},\tilde{u}_{i,t},N^v_{i,t}),$
where
\begin{equation}
\begin{aligned}
    \m{L}_{i,t}(x^k_{i,t},q_{ij,t},r_{ij,t},\lambda^k_{i,t},\tilde{u}_{i,t},N^v_{i,t})\Big|_{\Phi} = \ & \sigma^r_t \sum_{j\in\m{I}} p_{ij,t} q_{ij,t}  - (\sigma^c_t - \sigma^c_{t+1} + \eta_t) \sum_{k\in\m{K}} \gamma_k x^k_{i,t} - \sigma^r_t (\gamma_e + \gamma_p) \sum_{k\in\m{K}} \lambda^k_{i,t}\tilde{u}_{i,t} \\
    & + (\theta_t - \sigma^r_t\gamma_e) N^v_{i,t} + (\theta_t - \sigma^r_t\gamma_e)\sum_{j\in\m{I}} \left(q_{ij,t}w^p_{i,t} + q_{ij,t}t_{ij} + r_{ij,t}t_{ij} \right) \\
    & + \sum_{j\in\m{I}} (\mu_{i,t} - \mu_{j,t}) (q_{ij,t} + r_{ij,t}) - \theta_t \sum_{k\in\m{K}} \lambda^k_{i,t} \left[\dfrac{(1-\varepsilon)\Omega}{\omega} - \tilde{t}_i\right].
\end{aligned}
\end{equation}
Note that $\sigma^c_t = 0$ if $t > T$. Therefore, we can separately optimize a subproblem for zone $i$ and stage $t$:
\begin{subequations}\label{eqn:subproblem}
\begin{align}
    \t{maximize} \quad & \m{L}_{i,t}(x^k_{i,t},q_{ij,t},r_{ij,t},\lambda^k_{i,t},\tilde{u}_{i,t},N^v_{i,t})\Big|_{\Phi} \\
    \t{subject to} \quad 
    & N^v_{i,t} > \hat{N}^v, \\
    & q_{ij,t} > 0, \ \forall j\in\m{I}, \\
    & x^k_{i,0} \leq x^k_{i,t} \leq \hat{x}^k_{i}, \ \forall k\in\m{K}, \\
    & 0 \leq \lambda^k_{i,t} \leq \rho_k x^k_{i,t}, \ \forall k\in\m{K}, \\
    & \lambda^k_{i,t} (\tilde{c}^k_{i,t} - \tilde{u}_{i,t}) \leq \epsilon, \ \forall k\in\m{K}, \\
    & \tilde{c}^k_{i,t} - \tilde{u}_{i,t} \geq 0, \lambda^k_{i,t} \geq 0, \ \forall k\in\m{K}.
\end{align}
\end{subequations}
However, finding the global optimum for each subproblem remains challenging due to their nonconcave nature. Fortunately, we observe that each subproblem \eqref{eqn:subproblem} can be further divided into two smaller-scale subproblems that are easier to solve.

The first one is an optimization over $N^v_{i,t},q_{ij,t}$, and $r_{ij,t},\forall j\in\m{I}$:
\begin{equation}\label{eqn:subproblem 1st}
\begin{aligned}
    \underset{N^v_{i,t},q_{ij,t},r_{ij,t}}{\t{maximize}} \quad & \sigma^r_t \sum_{j\in\m{I}} p_{ij,t} q_{ij,t} + (\theta_t - \sigma^r_t\gamma_e)\sum_{j\in\m{I}} \left(q_{ij,t}w^p_{i,t} + q_{ij,t}t_{ij} + r_{ij,t}t_{ij} \right) \\
    & + (\theta_t - \sigma^r_t\gamma_e) N^v_{i,t} + \sum_{j\in\m{I}} (\mu_{i,t} - \mu_{j,t}) (q_{ij,t} + r_{ij,t}) \\
    \t{subject to} \quad & N^v_{i,t} \geq \hat{N}^v, \\
    & q_{ij,t} > 0,  r_{ij,t}\geq 0, \ \forall j\in\m{I}.
\end{aligned}
\end{equation}
Recall that the revenue term is concave with respect to passenger demand $q_{ij,t}$ for most commonly used demand functions (e.g, logit model, linear demand function, exponential demand function, etc). Therefore, given any fixed $N^v_{i,t}$, passenger waiting time is uniquely determined by \eqref{eqn:pickup time}, and the above problem becomes a concave program over $q_{ij,t}$ and $r_{ij,t}$, allowing us to find to the optimal values of $q_{ij,t}$ and $r_{ij,t}$ corresponding to different values of $N^v_{i,t}$. We can therefore optimally solve \eqref{eqn:subproblem 1st} by searching over a one-dimensional space of $N^v_{i,t}$.

The second subproblem is a nonconcave optimization over $\tilde{u}_{i,t}$ $x^k_{i,t}$ and $\lambda^k_{i,t},\forall k\in\m{K}$:
\begin{equation}\label{eqn:subproblem 2nd}
\begin{aligned}
    \underset{x^k_{i,t},\lambda^k_{i,t},\tilde{u}_{i,t}}{\t{maximize}} \quad & - (\sigma^c_t - \sigma^c_{t+1} + \eta_t) \sum_{k\in\m{K}} \gamma_k x^k_{i,t} - \sigma^r_t (\gamma_e + \gamma_p) \sum_{k\in\m{K}} \lambda^k_{i,t}\tilde{u}_{i,t}  - \theta_t \sum_{k\in\m{K}} \lambda^k_{i,t} \left[\dfrac{(1-\varepsilon)\Omega}{\omega} - \tilde{t}_i\right] \\
    \t{subject to} \quad &  x^k_{i,0} \leq x^k_{i,t} \leq \bar{x}^k_i, \ \forall k, \\
    & 0 \leq \lambda^k_{i,t} \leq \rho_k x^k_{i,t}, \ \forall k, \\
    & \lambda^k_{i,t} (\tilde{c}^k_{i,t} - \tilde{u}_{i,t}) \leq \epsilon, \ \forall k, \\
    & \tilde{c}^k_{i,t} - \tilde{u}_{i,t} \geq 0, \lambda^k_{i,t} \geq 0, \ \forall k.
\end{aligned}
\end{equation}
Given that the set $\m{K}$ is two-dimensional, this subproblem has a decision space of five dimensions. Despite its nonconcavity, we can proceed by adopting a hierarchical approach where we first fix $\tilde{u}_{i,t}$ and then solve the two residual problems over $(x^k_{i,t},\lambda^k_{i,t})$ for different $k$. Note that these two residual problems are independent and each one is two-dimensional with a closed and bounded decision space. Therefore, through the proposed hierarchical approach, the computational cost of (\ref{eqn:subproblem 2nd}) is equivalent to two three-dimensional optimization problems, which can be efficiently solved using brute-force enumeration without requiring the problem to be concave. Overall, we divide the subproblem into \eqref{eqn:subproblem 1st} and \eqref{eqn:subproblem 2nd}, each of which has favorable properties (concavity for \eqref{eqn:subproblem 1st} and low-dimensionality for \eqref{eqn:subproblem 2nd}) and can be solved using a hierarchical approach. Consequently, we can optimally solve each subproblem \eqref{eqn:subproblem}.

Based on the weak duality theorem, for any given Lagrange multipliers $\eta,\mu,\theta$, the solutions to all subproblems collectively provide a valid upper bound for the relaxed problem \eqref{eqn:relaxed problem}. This bound can be further improved by iteratively updating Lagrange multipliers within a dual decomposition framework as follows:
\begin{subnumcases}
    {\label{eqn:update multiplier}}
    \eta_t \leftarrow \max \left\{ 0, \eta_t - \zeta \left( b_t - \sum_{i\in\m{I}}\sum_{k\in\m{K}} \gamma_k x^k_{i,t} \right) \right\}, \\
    \mu_{i,t} \leftarrow \mu_{i,t} - \zeta \left( \sum_{j\in\m{I}} (q_{ij,t} + r_{ij,t}) - \sum_{j\in\m{I}} (q_{ji,t} + r_{ji,t}) \right) , \\
    \theta_{t} \leftarrow \theta_{t} - \zeta \left[ \sum_{i\in\m{I}} N^v_{i,t} + \sum_{i\in\m{I}}\sum_{j\in\m{I}} \left(q_{ij,t} w^p_{i,t} + q_{ij,t} t_{ij} + r_{ij,t} t_{ij}\right) - \sum_{i\in\m{I}}\sum_{k\in\m{K}} \lambda^k_{i,t} \left(\dfrac{(1-\varepsilon)\Omega}{\omega} - \tilde{t}_i\right) \right],
\end{subnumcases}
where $\zeta$ is the step size. Due to the model relaxation, the established upper bound is an overestimation of the globally optimal solution to the original problem \eqref{eqn:original problem}. In other words, the bound is also a theoretical upper bound for the original problem, which allows us to quantify the optimality of its solution. We summarize the derivation procedure as follows:
\begin{itemize}
    \item {\bf Step 1:} Solve the original problem \eqref{eqn:original problem} via nonlinear programming solvers. The solution provides a valid lower bound for \eqref{eqn:original problem}.
    \item {\bf Step 2:} Solve the relaxed problem \eqref{eqn:relaxed problem} via nonlinear programming solvers and obtain the corresponding Lagrange multipliers $(\eta,\mu,\theta)$.
    \item {\bf Step 3:} Plug $(\eta,\mu,\theta)$ into \eqref{eqn:subproblem}. Find the globally optimal solution to each subproblem by solving \eqref{eqn:subproblem 1st} and \eqref{eqn:subproblem 2nd}, respectively.
    \item {\bf Step 4:} Combine subproblems' solutions to establish a valid upper bound for \eqref{eqn:original problem}. Evaluate the bound performance by measuring the gap between lower and upper bounds.
    \item {\bf Step 5:} Terminate the procedure if the optimality gap is satisfactorily tight. Otherwise, update the multipliers according to \eqref{eqn:update multiplier} and go to {\bf Step 3}.
\end{itemize}
The aforementioned procedure can be terminated at any iteration, and regardless of when it is terminated, the derived bound is always a valid upper bound for the original problem \cite{bertsekas_nonlinear_2016}. The proposed algorithm is an example of how the problem's unique structure can be leveraged to assess the solution optimality of a highly nonconcave optimization problem after proper model relaxations, given that a feasible solution can be derived by some other gradient-based nonlinear programming solvers. We emphasize that the proposed method is nontrivial and rather useful. Because for most large-scale nonconcave problems, it may not be challenging to find a local solution using nonlinear programming solvers, but it is significantly difficult to examine how good the solution is compared to the unknown globally optimal one. In this work, we tackle this intrinsic challenge by relaxing the original problem to another nonconcave problem with decomposable structures. By leveraging the problem's favorable properties, we show that each subproblem can be optimally solved even in the presence of nonconcavity, thereby allowing us to evaluate the distance between the obtained solution and the globally optimal one.

\section{Numerical Studies}
Based on the proposed formulation and analytical analysis, this section presents several case studies to validate the solution optimality and discuss the effectiveness of the joint deployment scheme.

\subsection{Simulation Settings}
A series of case studies are conducted for Manhattan, New York City, where a platform jointly deploys charging infrastructure and operates a ride-hailing fleet. The island is divided into 20 zones.\footnote{The proposed model and solution framework are theoretically applicable to real-world networks with many more zones. But a 20-zone model should suffice to show the efficacy of our methodology.} Real data from the New York City Taxi and Limousine Commission (NYCTLC) \cite{tlc} was used to estimate travel demand $\bar{q}_{ij}$, including pickup and dropoff locations as well as trip duration. An exponential function is utilized as the demand model in \eqref{eqn:demand model}:
\begin{equation} \label{eqn:exponential model}
    q_{ij,t} = \bar{q}_{ij} \exp\left(-\alpha (p_{ij,t} + \beta w^p_{i,t})\right),
\end{equation}
where $\alpha > 0$ is a sensitivity parameter. We consider a three-stage infrastructure investment plan, i.e., $T = 3$. To capture the entire expansion process of an infrastructure network, we define $x_{i,0}^k =0, \forall i\in\m{I},k\in\m{K}$, implying that there is no charging facility before the first planning stage. Parameters related to EVs and infrastructure are calibrated to align with industry practices. Particularly, we assume all charging and battery swapping stations share the same configuration to ensure a fair comparison. For ease of comprehension, the number of swapping stations is used to quantify the total budget. On the other hand, parameters related to ride-hailing operations are calibrated via reverse engineering with partial reference to \cite{li_regulating_2019} such that the proposed model reproduces an equilibrium solution that is close to empirical data. The detailed parameter settings are summarized and justified in \ref{appendix:parameter}. In Steps 1 and 2 of our proposed solution procedure, we adopt the interior-point algorithm implemented in fmincon to obtain local solutions and Lagrange multipliers. Numerical simulations in this study are conducted on a desktop computer equipped with an Intel 8-core i7-9700 CPU.

\subsection{Solution Optimality}
\begin{figure}[!th]
    \centering
    \includegraphics[width=0.45\textwidth]{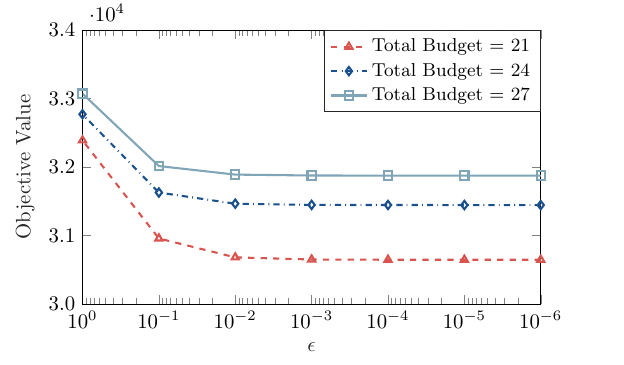}
    \caption{Objective value under different values of relaxation parameter $\epsilon$.}
    \label{fig:convergence}
\end{figure}
As mentioned in \ref{sec:stationary point}, a stationary point of the joint planning problem \eqref{eqn:original problem} can be found by iteratively reducing the relaxation parameter $\epsilon$. As shown in \autoref{fig:convergence}, we can obtain a local solution under different settings, which serves as a lower bound for \eqref{eqn:original problem}. While an upper bound can be established using the proposed relaxation and decomposition approach. The solution optimality is comprehensively evaluated across various budget levels, parameter settings, and network sizes. A six-zone network is included as it will be used in the subsequent section to visualize the expansion process of infrastructure networks. To benchmark the established upper bound, we derive another upper bound by applying Lagrangian relaxation on the original problem, where the complementarity constraint is linearized using the big-M method to enable decomposition. The derivation of this benchmark upper bound is detailed in \ref{appendix:benchmark}.

\begin{table}[ht!]
\fontsize{9}{12}\selectfont
\centering
\caption{Solution optimality under different simulation settings.}
\vspace{0.1cm}
\begin{tabular}{>{\centering\arraybackslash}m{1.1cm}>{\centering\arraybackslash}m{1.1cm}>{\centering\arraybackslash}m{1.1cm}>{\centering\arraybackslash}m{1.1cm}>{\centering\arraybackslash}m{1.3cm}>{\centering\arraybackslash}m{1.3cm}>{\centering\arraybackslash}m{1.3cm}>{\centering\arraybackslash}m{1.3cm}>{\centering\arraybackslash}m{1.3cm}>{\centering\arraybackslash}m{1.43cm}}
\toprule
\multirow{2}{*}{$\gamma$} & \multirow{2}{*}{$\gamma_p$} & Totl.  & \# of & LB & BM UB & \multirow{2}{*}{BM Gap} & Our UB & \multirow{2}{*}{Our Gap} & Totl.Time \\
& & Budget & Zones & $(\times 10^4)$   & $(\times 10^4)$ & & $(\times 10^4)$ & & (min) \\ \midrule
0.90 & 20 & 21 & 20 & 2.992 & 17.171 & 82.57\% & 3.088 & 3.11\% & 99.63  \\
0.90 & 20 & 24 & 20 & 3.045 & 16.965 & 82.05\% & 3.125 & 2.56\% & 104.12 \\
0.90 & 20 & 27 & 20 & 3.071 & 15.100 & 79.66\% & 3.156 & 2.71\% & 100.25 \\ \midrule\midrule
0.85 & 20 & 21 & 20 & 2.608 & 15.172 & 82.81\% & 2.693 & 3.15\% & 95.04 \\
0.85 & 20 & 24 & 20 & 2.652 & 15.262 & 82.62\% & 2.725 & 2.68\% & 96.03 \\
0.85 & 20 & 27 & 20 & 2.673 & 13.589 & 80.33\% & 2.755 & 2.97\% & 94.72 \\ \midrule\midrule
0.95 & 20 & 21 & 20 & 3.435 & 18.453 & 81.39\% & 3.543 & 3.07\% & 99.88 \\
0.95 & 20 & 24 & 20 & 3.498 & 18.129 & 80.70\% & 3.589 & 2.52\% & 99.79 \\
0.95 & 20 & 27 & 20 & 3.531 & 16.891 & 79.10\% & 3.620 & 2.48\% & 99.03 \\ \midrule\midrule
0.90 & 15 & 21 & 20 & 3.066 & 15.102 & 79.70\% & 3.164 & 3.09\% & 92.51 \\
0.90 & 15 & 24 & 20 & 3.117 & 14.615 & 78.67\% & 3.191 & 2.33\% & 99.58 \\
0.90 & 15 & 27 & 20 & 3.132 & 13.348 & 76.54\% & 3.218 & 2.69\% & 93.62 \\ \midrule\midrule
0.90 & 25 & 21 & 20 & 2.921 & 18.385 & 84.11\% & 3.021 & 3.31\% & 102.86 \\
0.90 & 25 & 24 & 20 & 2.976 & 18.233 & 83.68\% & 3.063 & 2.85\% & 101.97 \\
0.90 & 25 & 27 & 20 & 3.012 & 17.219 & 82.51\% & 3.097 & 2.74\% & 98.75 \\ \midrule\midrule
0.90 & 20 & 21 & 6  & 3.065 & 27.091 & 88.69\% & 3.193 & 4.03\% & 34.88 \\
0.90 & 20 & 24 & 6  & 3.145 & 29.496 & 89.34\% & 3.233 & 2.72\% & 35.72 \\
0.90 & 20 & 27 & 6  & 3.188 & 30.777 & 89.64\% & 3.263 & 2.33\% & 36.91 \\
\bottomrule
\end{tabular}
\label{table:solution optimality}
\end{table}

\autoref{table:solution optimality} reports the optimality gap and the computation time of Step 3 after one iteration, which is the most time-consuming component in the proposed method. The results demonstrate the consistent performance of our approach across different settings. Using the multipliers obtained in Step 2 without further updates, we can achieve a tight upper bound with an optimality gap of approximately 3\%, indicating that a near-optimal solution to \eqref{eqn:original problem} is attained. In contrast, the benchmark method cannot provide an accurate estimation of the unknown globally optimal solution. A major reason is that in the benchmark approach, the complementarity constraint has to be first linearized using the big-M method and then dualized to ensure decomposability. However, this process introduces strong linearity in the relaxed problem, which tends to produce extreme solutions that compromise the quality of the derived bound. As a result, the Lagrangian relaxation in the benchmark approach only yields a trivial upper bound with up to 80\% optimality gap, providing little insight into the problem's true complexity. Our proposed method, however, circumvents the need to dualize the complementarity constraint by employing appropriate relaxations based on the inherent characteristics of the problem. By doing so, we simplify the complex model into a more tractable form with favorable structures, while still capturing the essential features of the problem. Importantly, the complementarity constraint remains respected in the relaxed subproblem, allowing us to derive a significantly tighter upper bound compared to the benchmark approach. Besides, the average computation time to optimally solve each subproblem in the proposed method is less than 2 minutes. We emphasize that this solving time is satisfactorily short for an infrastructure planning problem that does not require real-time computation. In fact, the computation time is orders of magnitude shorter than the implementation time of a deployment plan (e.g., in the range of months or years), rendering it insignificant for the planning purpose. This comprehensive evaluation confirms the effectiveness and efficiency of our proposed approach.

\subsection{Synergy between Charging and Battery Swapping}
\begin{figure*}[!th]
    \centering
    \subfloat[Total Profit]{\includegraphics[width=0.35\textwidth]{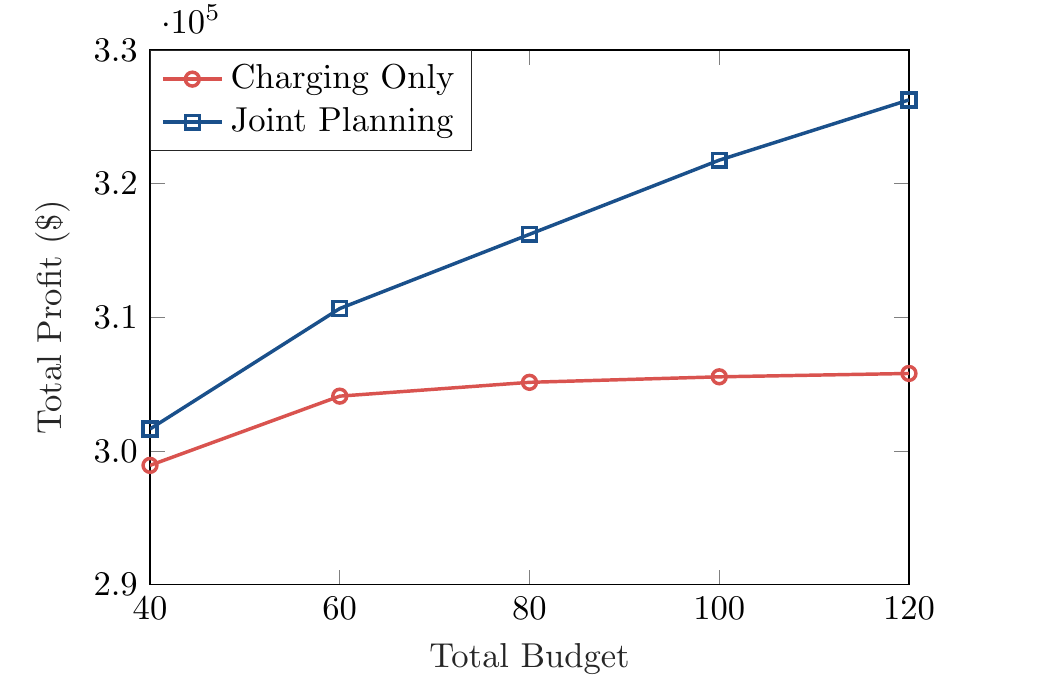}%
    \label{fig:total profit}}
    \hspace{-0.0cm}
    \subfloat[Average Charging Cost]{\includegraphics[width=0.35\textwidth]{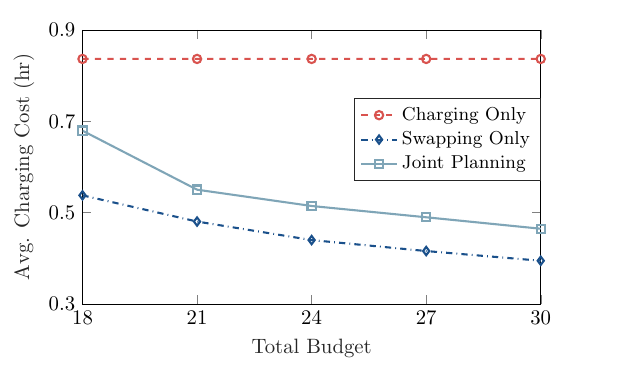}%
    \label{fig:average cost}}
    \caption{Performance of the infrastructure network under different deployment schemes.}
    \label{fig:synergy}
\end{figure*}

To assess the performance of joint planning, we vary the total budget and compare different deployment strategies, including jointly deploying charging and battery swapping facilities versus deploying only one of them. The total profit obtained from each strategy is shown in Figure \ref{fig:total profit}, and the corresponding average charging cost, defined as
\begin{equation}
    \bar{c} = \dfrac{\sum_{t\in\m{T}}\sum_{i\in\m{I}}\sum_{j\in\m{I}}\sum_{k\in\m{K}} f^k_{ij,t} u_{i,t}}{\sum_{t\in\m{T}}\sum_{i\in\m{I}}\sum_{j\in\m{I}}\sum_{k\in\m{K}} f^k_{ij,t} },
\end{equation}
is presented in Figure \ref{fig:average cost}. Additionally, we introduce the notion of cross-zone traffic to quantify the disparity between potential charging demand and actual charging demand. This metric, defined as
\begin{equation}
    d_t = \dfrac{1}{\sum_{i\in\m{I}} \Lambda_{i,t}} \left \| \Lambda_t - \sum\nolimits_{k\in\m{K}} \lambda^k_t  \right \|_2,
\end{equation}
provides a measure of the overall amount of relocation between zones due to drivers' charging behaviors.

\begin{figure*}[!th]
    \centering
    \subfloat[Total Budget = 21]{\includegraphics[width=0.35\textwidth]{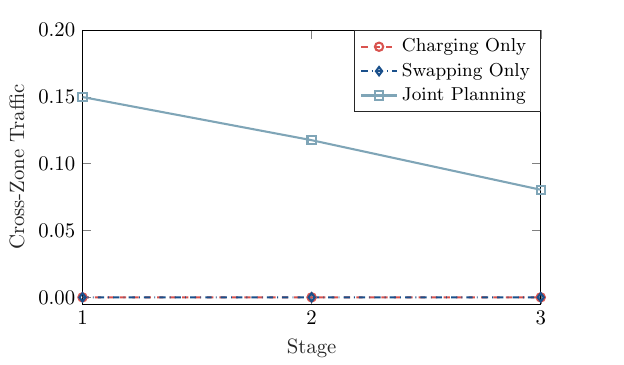}%
    \label{fig:crosszone traffic 7}}
    \hspace{-0.7cm}
    \subfloat[Total Budget = 24]{\includegraphics[width=0.35\textwidth]{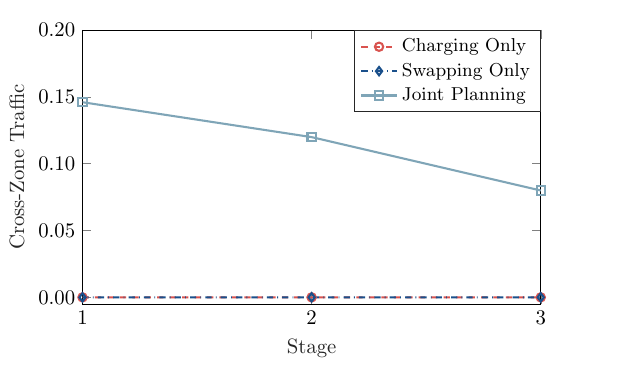}%
    \label{fig:crosszone traffic 8}}
    \hspace{-0.7cm}
    \subfloat[Total Budget = 27]{\includegraphics[width=0.35\textwidth]{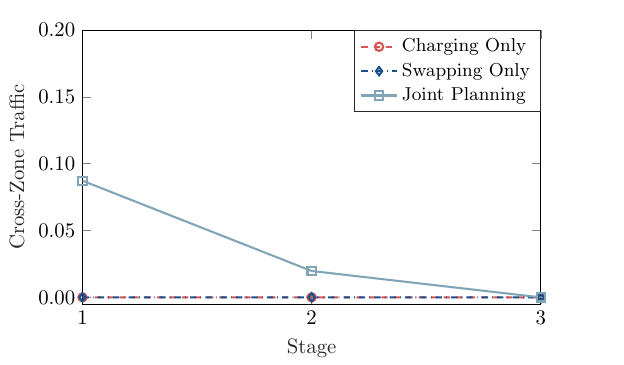}%
    \label{fig:crosszone traffic 9}}
    \caption{Distance between the potential charging demand and actual charging demand under different budget levels.}
    \label{fig:crosszone traffic}
\end{figure*}

On the one hand, battery swapping stations are costly and difficult to scale, while charging stations offer a more cost-effective solution for supporting a large EV fleet. Therefore, when the budget is limited, charging-only deployment yields a total profit that is 3.3\% higher than swapping-only deployment. On the other hand, the long charging time poses an inherent bottleneck on plug-in charging, while battery swapping naturally overcomes this bottleneck with its much quicker turnaround. As a result, with a sufficient budget that can support an extensive swapping network, swapping-only deployment can lead to a substantial profit improvement of 15.8\% and, at the same time, achieve a remarkable charging cost reduction of up to 52.8\% compared to charging-only deployment.

In contrast, joint planning combines the advantages of both charging and battery swapping stations, thereby harnessing the synergistic value between the two facilities. Under the joint planning scheme, the platform can not only address the scaling challenge of battery swapping by building charging stations at the early stage but also overcome the bottleneck of plug-in charging and enhance fleet utilization by building battery swapping stations. Consequently, joint planning outperforms other deployment strategies in terms of profit maximization. {However, the synergistic benefit of joint planning depends on the budget level as well as the scenario it is compared to. Specifically, under a tight budget, joint planning yields a total profit that is significantly higher than swapping-only deployment (i.e., by 11.7\%). In this case, the major bottleneck is the budget as it is too limited to afford an extensive swapping network and thus cannot support a large EV fleet. Conversely, when the budget is generous, the bottleneck becomes the low fleet utilization arising from time-consuming plug-in charging. The total profit of joint planning is significantly higher than charging-only deployment under a large budget (i.e., by 17.5\%). This is because by synergistically deploying charging and battery swapping stations, the platform can effectively reduce the charging cost and thus address this bottleneck.} As shown in Figure \ref{fig:average cost}, compared to charging-only deployment, joint planning yields a charging cost reduction of up to 44.4\% under a sufficient budget. As the budget increases, the charging cost of joint planning is even comparable to that of swapping only, which underscores the potential of joint planning to capitalize on the complementary strengths of charging and swapping facilities. {We formalize the superiority of joint planning as follows.}
\begin{proposition}
    {The total profit of joint planning is always larger than or  equal to that of charging-only or swapping-only deployment.}
\end{proposition}
\begin{proof}
    {Any charging-only deployment strategy can be interpreted as a joint planning strategy with the number of swapping stations being fixed at its initial value, i.e., $x^s_{i,t} = x^s_{i,0}, \forall t\in\m{T}$, which is always a feasible solution to the original optimization problem \eqref{eqn:original problem}. Therefore, the objective value under charging-only deployment will not be higher than that under joint planning. The same reasoning can be applied to swapping-only deployment as well.}
\end{proof}
{The underlying reason is that charging or swapping-only deployment can be considered as a special case of joint planning. Its performance is equal to that of joint planning only when the optimal joint planning strategy falls on the boundary, where only one type of facility is built.}

\begin{figure*}[!th]
    \centering
    \subfloat[Total Budget = 18]{\includegraphics[width=0.55\textwidth]{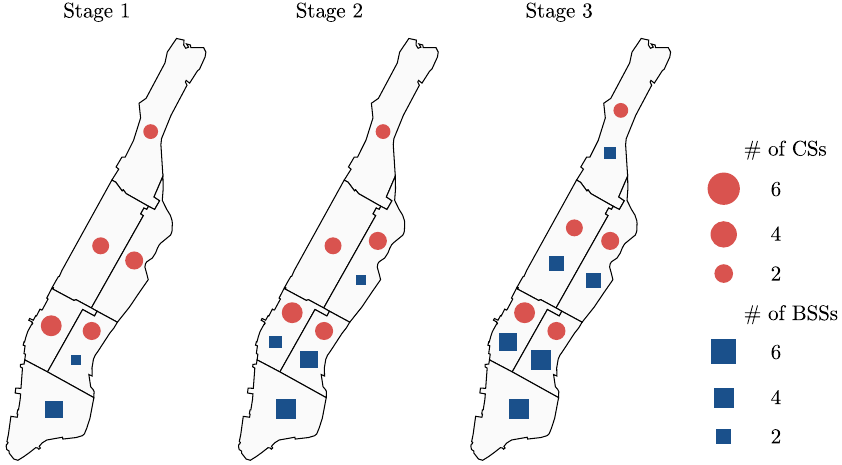}%
    \label{fig:expansion 6}}

    \subfloat[Total Budget = 24]{\includegraphics[width=0.55\textwidth]{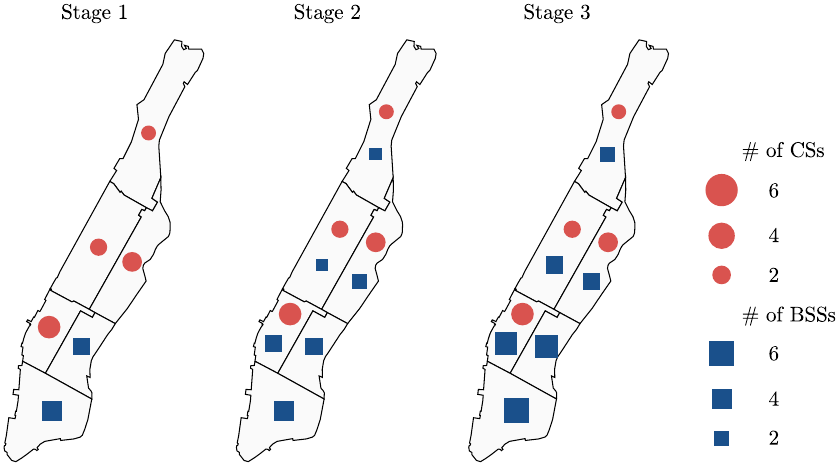}%
    \label{fig:expansion 8}}

    \subfloat[Total Budget = 30]{\includegraphics[width=0.55\textwidth]{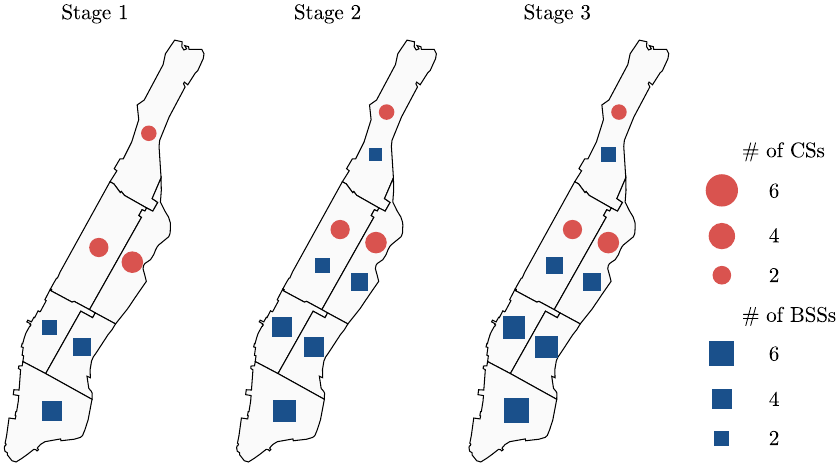}%
    \label{fig:expansion 10}}
    \caption{Infrastructure network expansion under different budget levels, with the red circles representing charging stations, the blue squares representing battery swapping stations, and their sizes representing the quantity of each type of facilities.}
    \label{fig:expansion}
\end{figure*}

In addition, we also observe that a multimodal charging network can result in a significant increase in cross-zone traffic compared to a unimodal charging network. The primary reason for this gap is the substantial difference between charging time and swapping time. Plug-in charging typically takes around an hour, while battery swapping can be completed within a few minutes. Under a unimodal charging network, drivers do not have strong incentives to travel long distances for recharging. However, given a multimodal charging network, drivers can reduce their charging costs by relocating to a neighboring zone with a battery swapping station as long as the resulting increase in travel time is shorter than the time gap between plug-in charging and battery swapping. Consequently, the volume of cross-zone traffic is higher under a joint planning strategy, particularly when the budget is insufficient to support an extensive battery swapping network, and  the volume gradually decreases as the budget increases. These findings highlight the potential impacts of infrastructure planning on urban traffic congestion and management, which should be carefully considered and addressed when deploying a charging infrastructure network.

To investigate the robustness of our findings, we perform a comprehensive sensitivity analysis, in which multiple key parameters, such as the cost ratio between the two facilities, the discount factor and the penalty parameter, are varied within a reasonable range. Detailed results of the sensitivity analysis, as presented in \ref{appendix:sensitivity analysis}, demonstrate the robustness of the superiority of joint planning. Across the various parameter settings tested, joint planning consistently outperforms alternative scenarios, which strengthens the argument for the effectiveness of joint planning.

\subsection{Infrastructure Network Expansion}

To further understand the benefits of joint planning, we present the expansion process of a mixed charging infrastructure network in \autoref{fig:expansion}. For better visualization, we re-partition Manhattan into 6 zones, with the upper three zones corresponding to areas with lower passenger demand and the lower three zones corresponding to areas with higher passenger demand.

It can be found that charging and battery swapping stations play different roles at different stages of long-term infrastructure planning. Due to the lower deployment cost, the platform prioritizes deploying charging stations at the early stage such that a dense charging network can be built to support a large EV fleet. As investment increases, the bottleneck of plug-in charging becomes increasingly apparent, necessitating an upgrade of the infrastructure network. The platform then deploys battery swapping stations to enhance fleet utilization. The expansion process highlights the complementary effect between the two facilities. 

Another important observation is that the deployment priority of battery swapping stations is always assigned to areas with higher passenger demand. This is because charging demand is positively related to passenger demand in each zone. High-demand areas typically have high charging demand. In this light, the platform prefers to first build swapping facilities in high-demand areas where more EVs could use battery swapping for recharge, so the economical benefit associated with quicker turnaround and improved fleet utilization can be fully harnessed. { Furthermore, by comparing Figure \ref{fig:expansion 6}-\ref{fig:expansion 10}, we can observe that the larger the budget level is, the fewer charging stations that will be built in the high-demand areas. These results also suggest that the platform should be far-sighted and avoid over-building charging stations in high-demand areas during the early stage, as this may lead to regret in the future when the budget is more abundant or when battery swapping is introduced.} Overall, these discussions shed light on how demand and budget levels will affect the optimal deployment strategy and thus can advise infrastructure planning.

{
\section{Extensions}
\subsection{Spatial Charging Equilibrium with Rejection}
Since both the charging and battery swapping stations have finite capacities, newly arriving EVs may possibly get rejected if the queue is full. The proposed mathematical model is based on simplified settings where the queue capacity is assumed to be sufficiently large, resulting in a negligible blocking probability. However, in order to provide a more comprehensive model, we consider the scenario where rejection is taken into account. Let $\upsilon^k_{i,t}$ represent the blocking probability for mode $k$. All the rejected EVs should be relocated by choosing a charging mode and location again. The flow conservation conditions of these rejected EVs can be represented as:
\begin{subnumcases}{}
    \sum_{k\in\m{K}} \upsilon^k_{i,t} \lambda^k_{i,t} = \sum_{j\in\m{I}} \sum_{k\in\m{K}} \hat{f}^k_{ij,t}, \\
    \hat{\lambda}^k_{i,t} = \sum_{j\in\m{I}} \hat{f}^k_{ji,t},
\end{subnumcases}
where $\hat{f}^k_{ij,t}$ is the flow of the relocated EVs that are rejected due to limited queuing capacity, and $\hat{\lambda}^k_{ij,t}$ represents the arrival rates of relocated rejected flows to zone $i$. Therefore, the average arrival rate at each station is given by
\begin{equation}
    \bar{\lambda}^k_{i,t} = \dfrac{1}{x^k_{i,t}} \left(\lambda^k_{i,t} + \hat{\lambda}^k_{i,t}\right),
\end{equation}
which determines the waiting time $w^k_{i,t}$ and the blocking rate $\upsilon^k_{i,t}$ according to the queueing model. In this case, the charging cost is then defined as
\begin{equation}
    c^k_{ij,t} = t_{ij} + (1-\upsilon^k_{j,t})\left(t_k + w^k_{j,t}\right) + \upsilon^k_{j,t} u_{j,t},
\end{equation}
where $u_{i,t}$ denotes the equilibrium charging cost for EVs from zone $i$. The spatial charging equilibrium is expressed by the following set of equations:
\begin{subnumcases}{}
    f^k_{ij,t},\hat{f}^k_{ij,t} \geq 0, \\
    c^k_{ij,t} - u_{i,t} \geq 0, \\
    f^k_{ij,t} (c^k_{ij,t} - u_{i,t}) = 0, \\
    \hat{f}^k_{ij,t} (c^k_{ij,t} - u_{i,t}) = 0.
\end{subnumcases}
Note that the model implicitly assumes that each charging demand will not be blocked for more than once. However, in practice, an EV may experience multiple rejections before he/she successfully joins a queue (although the chance is rather small). One can easily account for the case of repeated rejections by introducing additional $\hat{f}^k_{ij,t}$ and $\hat{\lambda}^k_{i,t}$ following the same definitions. This model serves as a ``first-order'' approximation of the true case, which is sufficiently accurate given the observation that EVs will proactively choose a facility with a short waiting time and, consequently, a low blocking rate. The extended model is validated in a numerical study and the results are presented in \ref{appendix:extension}. We find that (1) the first-order approximation adequately captures the charging rejection, with minimal discrepancies when compared to higher-order approximations; and (2) the simplified model aligns closely with the first-order approximation, indicating that overlooking charging rejection does not alter the qualitative outcomes of our analysis.

\subsection{Spatial Charging Equilibrium with Differentiated Pricing}
Considering the significantly shorter service time of battery swapping, EV drivers tend to over-utilize swapping stations, leading to extended waiting times such that the charging costs of the two modes are the same. It is reasonable to implement differentiated service fees for charging and battery swapping stations as a means to regulate demand. In light of this, we extend our proposed model to incorporate a differentiated pricing strategy. Let $s^k_{i,t}$ denote the service fee (measured in time) for mode $k$ in zone $i$ at stage $t$. Consequently, the charging cost is re-defined as
\begin{equation}
    c^k_{i,t} = t_{ij} + t_k + w^k_{j,t} + s^k_{j,t}.
\end{equation}
By directly adjusting the service fees for different modes across various locations, the platform can induce the desired charging behaviors. Our primary focus is to examine how the price differentiation between battery swapping and plug-in charging will impact market outcomes. To facilitate a fair comparison between the simplified model and the extended model, we normalize the charging price to zero and use the price gap between battery swapping and plug-in charging, denoted as $s_{i,t} \geq 0$, as a decision. The positivity of $s_{i,t}$ arises from the observation that battery swapping is typically more expensive than plug-in charging. As such, the charging cost is modified as follows:
\begin{subnumcases}
    {c^k_{i,t} = }
    t_{ij} + t_k + w^k_{j,t}, & $\t{if } k = c$, \\
    t_{ij} + t_k + w^k_{j,t} + s_{j,t}, & $\t{if } k = s$.
\end{subnumcases}
If $s_{i,t} = 0$, the model degrades to the original one where EV drivers' decisions only depend on time costs. However, the platform may also set a positive price gap, i.e., $s_{i,t} > 0$, if it finds that the benefits of demand regulation outweigh the penalties associated with the increased equilibrium charging cost. We validate this extended model and compare its performance to the case without differentiated pricing in \ref{appendix:extension}. Our results indicate that the implementation of differentiated pricing enhances the platform’s flexibility in charging demand management, thereby leading to an improved profit. However, the improvement is insignificant, which is less than 1\% for a large range of total budgets.
}

\section{Conclusion}
This study investigates the planning of a multimodal charging network, where charging stations and battery swapping stations are jointly deployed to support an electric ride-hailing fleet. A multi-stage charging network expansion model is formulated to capture how the platform makes infrastructure planning and operational decisions to maximize it profit. The model incorporates demand elasticity, spatial charging equilibrium, charging and swapping congestion, and other fundamental components in an electric ride-hailing market. A relaxed formulation is developed to establish a theoretical upper bound for the nonconcave joint deployment problem. Moreover, extensive numerical studies are conducted to validate the proposed model using real data from Manhattan.

Our findings indicate that joint planning charging and battery swapping stations can elicit synergistic value between the two facilities and yield a "one plus one more than two" benefit. Compared to deploying only one of them, joint planning not only addresses the scaling issue of battery swapping but also mitigates the inherent bottleneck associated with plug-in charging. At the early stage of infrastructure deployment, the platform prioritizes a dense network of charging stations that can support a large EV fleet. As the investment budget increases, the platform starts building battery swapping stations to improve fleet utilization and maximize profit. These results confirm the effectiveness of the proposed multimodal charging network. However, a multimodal charging network will lead to a larger volume of cross-zone traffic than a unimodal charging network, which highlights the potential impacts of infrastructure planning on traffic congestion and management. We also briefly discuss the deployment priority of battery swapping stations and show that this priority always assigned to high-demand area where more EVs can adopt battery swapping for recharging.

This study opens up several research directions that await further investigation. For instance, we focus on the steady-state performance in this work, whereas passenger demand varies over time in practice, and charging demand also changes accordingly. Therefore, a possible extension is to jointly consider the long-run planning decisions and real-time operational decisions. This can account for the temporal dynamics of passenger demand and charging demand, leading to more robust and efficient infrastructure planning solutions. Besides, drivers decide when to charge depending on the current state of charge, time-varying passenger demand, and the availability of charging facilities. It is therefore intriguing to explore how a mixed infrastructure network will reshape drivers' charging schedule when considering these real-time factors, {and how to extend the model to explicitly characterize the interdependence between stationary charging demand and drivers' charging behaviors, which is a limitation of the current model.} In addition, this paper considers a ride-hailing platform collaborates with infrastructure providers to deploy charging infrastructure and demonstrates the economical benefits of such a partnership. While it remains unclear how these benefits would be divided between the two participants and how to ensure the stability of this partnership. These follow-up questions warrant further study. Moreover, future electric ride-hailing fleets could pose challenges and opportunities to urban power grids, especially considering the energy storage feature of battery swapping stations. So a dynamic charging and dispatching policy with integration of renewable energies is another interesting extension.

\section*{Acknowledgments}  {This research was supported by Hong Kong Research Grants Council under project 26200420, 16202922, and National Science Foundation of China under project 72201225.}

\bibliographystyle{unsrt}
\bibliography{reference}

\appendix

\section{Parameter Settings}
\label{appendix:parameter}
To ensure a realistic representation of charging infrastructure, we calibrate the relevant parameters based on industry practices. For instance, it is suggested that around 70\% of charging stations in US metropolitan areas have less than six chargers \cite{nicholas_estimating}. Therefore, we assume each charging station in our model has five chargers. Furthermore, taking NIO's battery swapping station as a reference, we set the number of swapping slots to 1 \cite{scmp_chinese}. In particular, charging and swapping stations are assumed to have the same configuration for a fair comparison. We thus consider each battery swapping station has five chargers and five batteries in circulation. The average charging time is 40 minutes, while each swapping session lasts 5 minutes. To minimize the impact of blocking, we assume both facilities have a queue capacity of 50. The upper bounds for the arrival rate at both facilities, i.e., $\rho_k$, are selected to result in an average waiting time of 1 hour (See \ref{appendix:queueing model} for more details). {In terms of deployment costs, we decide the value based on infrastructure configuration. An estimation from 2019 \cite{nicholas_estimating} indicates that the cost of deploying each 50kW DC fast charger in U.S. metropolitan areas is approximately \$40,000, depending on the location, technology used, and other factors. Hence, the total cost of a charging station with five chargers amounts to \$200,000 (note that the cost of the charging station depends on its configuration, i.e., number of chargers). We apportion this cost to an hourly basis over a one-year period and assume $\gamma_c = \$20$.  While building each battery swapping station costs around \$500,000 according to NIO \cite{nedelea_take}. Taking into account NIO's technological advantages and economies of scale, as well as higher land acquisition costs in Manhattan, we assume the cost of building a battery swapping station is five times that of a charging station with the same configuration, i.e., five chargers and five batteries, so that the cost of battery swapping station also depends on its configuration (i.e., number of chargers inside).} The operational cost of EVs is set to be \$25/hr based on the hourly earnings of Uber drivers in NYC \cite{cnbc_make}. We consider each stage spans one year and set the discount factor as 0.9. Additionally, instead of tuning $(\varepsilon,\omega,\Omega)$ separately in \eqref{eqn:energy balance}, we fix $(1-\varepsilon)\Omega/\omega = 8$, meaning that EVs can sustain for eight hours between two top-ups. This is close to the charging pattern of a large-scale electric taxi fleet \cite{lei_understanding_2022}.

\begin{table}[h!]
\fontsize{9}{12}\selectfont
\centering
\caption{Summary of model parameters}
\vspace{0.1cm}
\begin{tabular}{llc}
\toprule
Parameter & Description & Value \\ \midrule
$\gamma_e$ & operational cost of each EV & \$25/hr \\
$\gamma_c$ & deployment cost of each charging station & \$20/hr \\
$\gamma_s$ & deployment cost of each battery swapping station & \$100/hr \\
$\gamma_p$ & penalty parameter for the total charging cost & \$20/hr \\
$\gamma$ & discount factor in the objective function & 0.9 \\
$\phi$ & parameters in the pickup time model & 1/3 \\
$\alpha$ & sensitivity parameter in the demand model & 0.12 \\
$\beta$ & passengers' value of time & \$90/hr \\
$t_c$ & average charging time & 40 min \\
$t_s$ & constant swapping time & 5 min  \\
$C$ & number of chargers in each station & 5 \\
$S$ & number of swapping slots in each station & 1 \\
$B$ & number of batteries at each swapping station & 5 \\
$K$ & queue capacity of each station & 50 \\
$\rho_c$ & upper bound for the arrival rate at charging stations & 6.767 \\
$\rho_s$ & upper bound for the arrival rate at swapping stations & 5.499 \\
$T$ & number of stages in the planning horizon & 3 \\
$H$ & number of stages in the operation horizon & 6 \\
$(1-\varepsilon)\Omega/\omega$ & charging frequency of each EV & 8 \\ \bottomrule
\end{tabular}
\label{table:parameter}
\end{table}

\section{Stationary Distribution of Swapping Queue}
\label{appendix:swapping queue}
As indicated in \cite{tan_queueing_2014}, the steady-state distribution of the swapping queue can be conveniently computed in matrix form without solving a linear equation system. Specifically, we can obtain $G(i,j)$ by simply solving $G \Pi = G$, where $G = [G(0,0),\cdots,G(0,B),\cdots,G(K,0),\cdots,G(K,B)]$ is a row vector and $\Pi$ is the transition matrix of the embedded Markov chain. To construct the matrix $\Pi$, we define the following $(B+1)\times(B+1)$ matrix $\Delta$:
\begin{equation}
    \Delta = \left[\begin{array}{ccc:ccc}
    f(0,0) & \cdots & 0 & 0 & \cdots & 0 \\
    \vdots & \ddots & \vdots & \vdots & \vdots & \vdots \\
    f(C,0) & \cdots & f(C,C) & 0 & \cdots & 0 \\ \hdashline
    f(C,0) & \cdots & f(C,C) & 0 & \cdots & 0 \\
    \vdots & \ddots & \vdots & \vdots & \vdots & \vdots \\
    f(C,0) & \cdots & f(C,C) & 0 & \cdots & 0 \\
    \end{array}\right],
\end{equation}
where $f(i,j)$ is given by \eqref{eqn:battery circulation}. As such, any infeasible transition can be easily excluded by 0 entries in the matrix. Then the transition matrix $\Pi$ can be constructed as follows:
\begin{algorithm}[ht!]
\fontsize{9}{12}\selectfont
\SetAlgoLined
\caption{Construction of matrix $\Pi$}
\label{alg:construc Pi}
\For{$i \leftarrow 1 \t{ to } (K+1)\times(B+1)$}{
\For{$j \leftarrow 1 \t{ to } (K+1)\times(B+1)$}{
$a = \left\lfloor (i-1)/(B+1) \right\rfloor$ \\
$b = (i-1) \mod (M+1)$ \\
$m = \left\lfloor (j-1)/(B+1) \right\rfloor$ \\
$n = (j-1) \mod (M+1)$ \\
$\delta=\begin{cases}
    g\left(m-1-a+\min\{a,b,S\}\right), & \t{if $m<K-1$}\\
    1 - \sum_{l=0}^{K-1} g\left(l-a+\min\{a,b,S\}\right), & \t{if $m=K-1$}
\end{cases}$ \\
$\Pi(i,j) = \delta \cdot \Delta(B-b+1,n-b+\min\{a,b,S\})$
}}
\end{algorithm}

\section{Performance of Queueing Models}
\begin{figure*}[!th]
    \centering
    \subfloat[Waiting Time]{\includegraphics[width=0.35\textwidth]{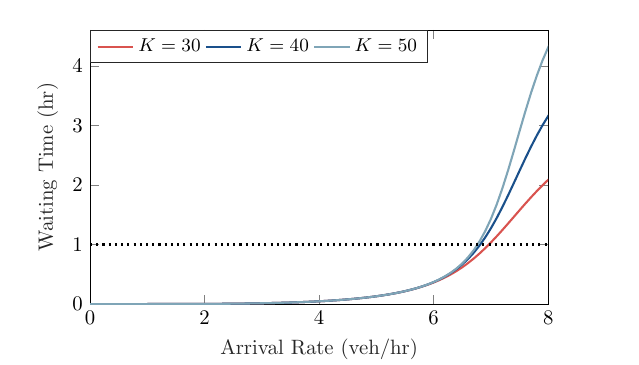}%
    \label{fig:charging waiting time}}
    \hspace{-0.0cm}
    \subfloat[Blocking Rate]{\includegraphics[width=0.35\textwidth]{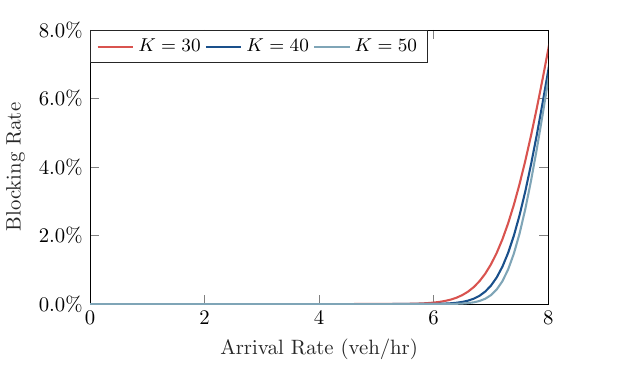}%
    \label{fig:charging blocking rate}}
    \caption{Queueing performance of a charging station.}
    \label{fig:charging queue}
\end{figure*}

\begin{figure*}[!th]
    \centering
    \subfloat[Waiting Time]{\includegraphics[width=0.35\textwidth]{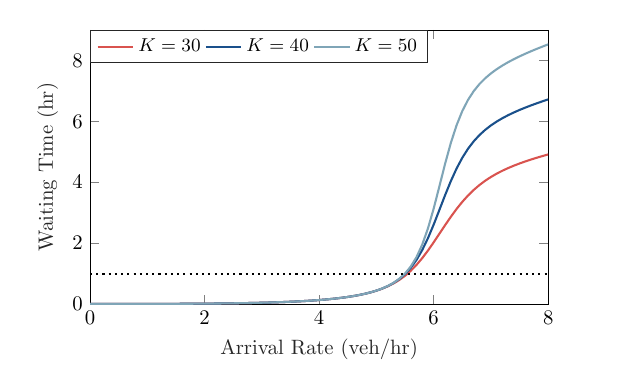}%
    \label{fig:swapping waiting time}}
    \hspace{-0.0cm}
    \subfloat[Blocking Rate]{\includegraphics[width=0.35\textwidth]{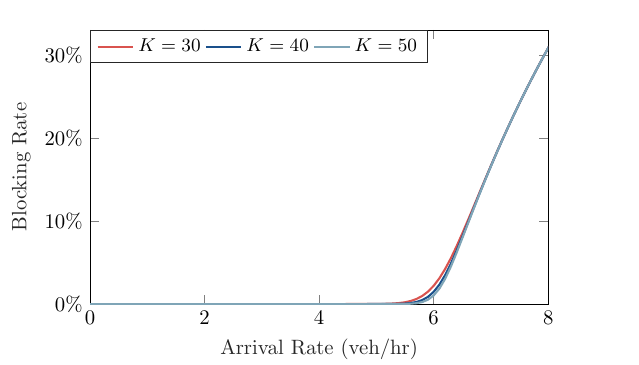}%
    \label{fig:swapping blocking rate}}
    \caption{Queueing performance of a battery swapping station.}
    \label{fig:swapping queue}
\end{figure*}

\autoref{fig:charging queue} and \autoref{fig:swapping queue} show the waiting time and blocking rate of both charging and swapping queues under different EV arrival rates. It can be found that when the waiting time is in a tolerable range, e.g., less than 1 hour, the blocking rate is negligible. To ensure acceptable waiting times, we impose constraint \eqref{eqn:demand upper bound}, which sets an upper bound on the EV arrival rate. This constraint limits the arrival rate such that the resulting waiting time is kept below 1 hour. As summarized in \autoref{table:blocking rate}, with a sufficiently large queue capacity, the blocking rate is marginal, e.g., less than 0.1\%, for any arrival rate below its upper bound. These results illustrate that we can effectively prevent the case of a full queue by properly limiting the maximum EV arrival rate and meanwhile choosing a suitable queue capacity.

\label{appendix:queueing model}

\begin{table}[ht!]
\fontsize{9}{12}\selectfont
\centering
\caption{EV arrival rate and the blocking rate when the waiting time is 1 hour.}
\vspace{0.1cm}
\begin{tabular}{|c|c|c||c|c|}
\hline
$K$  & $\rho_c$  & Blocking Rate      & $\rho_s$  & Blocking Rate      \\ \hline
30 & 6.965 & 1.055\% & 5.546 & 0.278\% \\
40 & 6.814 & 0.255\% & 5.510 & 0.063\% \\
50 & 6.767 & 0.074\% & 5.499 & 0.016\% \\ \hline
\end{tabular}
\label{table:blocking rate}
\end{table}

\section{Derivation of Benchmark Upper Bounds}
\label{appendix:benchmark}
Aside from the proposed method, another upper bound is established as a benchmark via Lagrangian relaxation. To obtain a decomposable structure from the original problem, we first insert \eqref{eqn:scaling n} into \eqref{eqn:demand distribution}:
\begin{equation}\label{eqn:rewritten transition matrix}
    \sum_{i\in\m{I}} \Lambda_{i,t} P_{ij,t} = \Lambda_{j,t},
\end{equation}
which is clearly separable over $i$ because $P_{ij,t}$ only depends on $N^v_{i,t},q_{ij,t},r_{ij,t},\forall j\in\m{I}$. Note that the energy balance condition also has a separable structure. For ease of notation, we combine \eqref{eqn:rewritten transition matrix} and \eqref{eqn:energy balance} as $H_{v,t} = \sum_{i\in\m{I}} H_{vi,t} = 0$, where $H_{vi,t}$ is a function of $\{N^v_{i,t},q_{ij,t},r_{ij,t},f^k_{ij,t}\}_{\forall j\in\m{I},k\in\m{K}}$. Moreover, we linearize the complementarity constraint \eqref{eqn:complementarity} by introducing a binary variable $z^k_{ij,t}$:
\begin{subnumcases}{\label{eqn:linear complementarity}}
    z^k_{ij,t} \in \{0,1\}, \\
    f^k_{ij,t} \leq M(1-z^k_{ij,t}), \\
    c^k_{ij,t} - u_{i,t} \leq M z^k_{ij,t}, \label{eqn:big M complementarity}
\end{subnumcases}
where $M$ is a sufficiently large constant. This is equivalent to the original complementarity constraint in the sense that when $f^k_{ij,t} > 0$, $z^k_{ij,t}$ must be 0, and thus $c^k_{ij,t} = u_{i,t}$. The equilibrium conditions \eqref{eqn:equilibrium condition} then becomes separable since $c^k_{ij,t}$ only depends on $\lambda^k_{j,t}$. To ensure the satisfaction of constraint qualification, we further relax \eqref{eqn:big M complementarity} as in \eqref{eqn:relaxed complementarity}:
\begin{equation}
    c^k_{ij,t} - u_{i,t} \leq M z^k_{ij,t} + \epsilon.
\end{equation}
As such, the original optimization problem becomes a mixed-integer nonlinear program that is subject to the following coupled constraints:
\begin{subequations}
\begin{align}
    & \sum\nolimits_{i\in\m{I}}\sum\nolimits_{k\in\m{K}} \gamma_k x^k_{i,t} \leq b_t, \ \forall t\in\m{T}, \label{eqn:con1} \\
    & x^k_{i,t} - x^k_{i,t-1} \geq 0, \ \forall t\in\m{T},i\in\m{I},k\in\m{K}, \label{eqn:con2} \\
    & \lambda^k_{i,t} = \sum\nolimits_{j\in\m{I}} f^k_{ji,t}, \ \forall t\in\m{T},i\in\m{I},k\in\m{K}, \label{eqn:con3} \\
    & H_{v,t} = \sum\nolimits_{i\in\m{I}} H_{vi,t} = 0, \ \forall t\in\m{T},v\in\m{V}, \label{eqn:con4} \\
    & c^k_{ij,t} - u_{i,t} \geq 0, \ \forall t\in\m{T},i\in\m{I},j\in\m{I},k\in\m{K}, \\
    & c^k_{ij,t} - u_{i,t} \leq M z^k_{ij,t} + \epsilon, \ \forall t\in\m{T},i\in\m{I},j\in\m{I},k\in\m{K}, \\
    & \sum\nolimits_{j\in\m{I}} (q_{ij,t} + r_{ij,t}) = \sum\nolimits_{j\in\m{I}} (q_{ji,t} + r_{ji,t}), \ \forall t\in\m{T},i\in\m{I}. \label{eqn:con5}    
\end{align}
\end{subequations}
Let $\Phi = (\eta,\delta,\nu,\theta,\psi,\phi,\mu)$ denote the Lagrange multipliers associated with these constraints. We then have the separable partial Lagrangian as follows:
\begin{equation}
\begin{aligned}
    \m{L} =& \sum_{t\in\m{T}} \left(\sigma^r_t \Pi^r_t - \sigma^c_t \Pi^c_t \right) - \sum_{t\in\m{T}} \sum_{i\in\m{I}}\sum_{k\in\m{K}} \eta_t\gamma_k x^k_{i,t} + \sum_{t\in\m{T}}\sum_{i\in\m{I}}\sum_{k\in\m{K}} \nu^k_{i,t} \Big(\lambda^k_{i,t}  - \sum_{j\in\m{I}} f^k_{ji,t}\Big) \\
    & + \sum_{t\in\m{T}}\sum_{i\in\m{I}}\sum_{j\in\m{I}}\sum_{k\in\m{K}} \psi^k_{ij,t} \Big(c^k_{ij,t} - u_{i,t}\Big) - \sum_{t\in\m{T}}\sum_{i\in\m{I}}\sum_{j\in\m{I}}\sum_{k\in\m{K}} \phi^k_{ij,t} \Big(c^k_{ij,t} - u_{i,t} - M z^k_{ij,t} \Big) \\
    & + \sum_{t\in\m{T}}\sum_{i\in\m{I}}\sum_{k\in\m{K}} \delta^k_{i,t} (x^k_{i,t} - x^k_{i,t-1}) + \sum_{t\in\m{T}}\sum_{i\in\m{I}}\sum_{v\in\m{V}} \theta_{v,t} H_{vi,t} + \sum_{t\in\m{T}}\sum_{i\in\m{I}} \sum_{j\in\m{I}} (\mu_{i,t} - \mu_{j,t}) (q_{ij,t} + r_{ij,t})
\end{aligned}
\end{equation}
The decomposable structure of this Lagrangian allows us to define multiple subproblems, each of which is a function of $(x^k_{i,t},u_{i,t},\Lambda_{i,t},\lambda^k_{i,t},N^v_{i,t},q_{ij,t},r_{ij,t},f^k_{ij,t},z^k_{ij,t}),\forall j\in\m{I},k\in\m{K}$. Furthermore, each subproblem can be divided into two classes of smaller-scale problems that can be independently solved. 

The first class is a bi-variate optimization of $x^k_{i,t}$ and $\lambda^k_{i,t}$:
\begin{equation}\label{eqn:subproblem1 benchmark}
\begin{aligned}
    \underset{x^k_{i,t},\lambda^k_{i,t}}{\t{maximize}} \quad & \left[\gamma_k(\sigma^c_{t+1} - \sigma^c_t - \eta_t) + (\delta^k_{i,t} - \delta^k_{i,t+1})\right] x^k_{i,t} + \nu^k_{i,t} \lambda^k_{i,t} \\
    & - \sigma^r_t \gamma_e \lambda^k_{i,t}(w^k_{i,t} + t_k) + \sum_{j\in\m{I}} (\psi^k_{ji,t} - \phi^k_{ji,t}) w^k_{i,t}\\
    \t{subject to} \quad & x^k_{i,0} \leq x^k_{i,t} \leq \bar{x}^k_i, \\
    & 0 \leq \lambda^k_{i,t} \leq \rho_k x^k_{i,t}.    
\end{aligned}
\end{equation}
Such a two-dimensional maximization problem can be efficiently solved even in the absence of concavity.

The second class is a mixed-integer nonlinear program over $(u_{i,t},\Lambda_{i,t},N^v_{i,t},q_{ij,t},r_{ij,t},f^k_{ij,t},z^k_{ij,t}),\forall j\in\m{I},k\in\m{K}$:
\begin{equation}\label{eqn:subproblem2 benchmark}
\begin{aligned}
    \underset{\substack{u_{i,t},\Lambda_{i,t},N^v_{i,t},q_{ij,t},\\ r_{ij,t},f^k_{ij,t},z^k_{ij,t}}}{\t{maximize}} \quad & \sigma^r_t \sum_{j\in\m{I}} p_{ij,t} q_{ij,t} -\sigma^r_t \gamma_e \sum_{j\in\m{I}} q_{ij,t}w^p_{i,t} -\sigma^r_t \gamma_e \sum_{j\in\m{I}}\Big(  q_{ij,t} + r_{ij,t} + \sum_{k\in\m{K}} f^k_{ij,t} \Big) t_{ij}  \\
    & -\sigma^r_t\gamma_e N^v_{i,t} -\sigma^r_t\gamma_p \Lambda_{i,t}u_{i,t} + \sum_{v\in\m{V}} \theta_{v,t} H_{vi,t} + \sum_{j\in\m{I}} (\mu_{i,t} - \mu_{j,t}) (q_{ij,t} + r_{ij,t}) \\
    & - \sum_{j\in\m{I}}\sum_{k\in\m{K}} \nu^k_{j,t} f^k_{ij,t}  + \sum_{j\in\m{I}}\sum_{k\in\m{K}} (\phi^k_{ij,t} - \psi^k_{ij,t}) u_{i,t} + M\sum_{j\in\m{I}}\sum_{k\in\m{K}} \phi^k_{ij,t} z^k_{ij,t} \\
    \t{subject to} \quad & N^v_{i,t} > \hat{N}^v, \\
    & q_{ij,t} > 0, \ \forall j\in\m{I}, \\
    & \Lambda_{i,t} = \sum\nolimits_{k\in\m{K}}\sum\nolimits_{j\in\m{I}} f^k_{ij,t}, \\
    & z^k_{ij,t} \in \{0,1\}, \ \forall j\in\m{I},k\in\m{K}, \\
    & f^k_{ij,t} \leq M(1-z^k_{ij,t}), \ \forall j\in\m{I},k\in\m{K}.
\end{aligned}
\end{equation}
We notice that complexities of the objective function mainly arise from $H_{vi,t}$. However, if $\Lambda_{i,t}$, $N^v_{i,t}$, and $\sum_{j\in\m{I}}q_{ij,t}w^p_{i,t}+\sum_{j\in\m{I}} (q_{ij,t} + r_{ij,t}) t_{ij}$ are fixed, $H_{vi,t}$ is simplified to a linear function of $q_{ij,t}$, $r_{ij,t}$ and $f^k_{ij,t}$. This problem consequently becomes a concave program over $q_{ij,t}$ and a mixed-integer linear program over other decisions. We can therefore perform a brute-force search over a three-dimensional space and then efficiently solve the two residual programs using off-the-shelf solvers.

Following the principle of Lagrangian relaxation, we obtain a valid upper bound for the original problem by optimally solving \eqref{eqn:subproblem1 benchmark} and \eqref{eqn:subproblem2 benchmark} for all $i\in\m{I}$ and $t\in\m{T}$. This bound is utilized as a benchmark to evaluate the performance of our proposed bound.

\newpage
\section{Results of Sensitivity Analysis}
\label{appendix:sensitivity analysis}
\begin{figure*}[!th]
    \centering
    \subfloat[$\gamma_s/\gamma_c = 3$]{\includegraphics[width=0.35\textwidth]{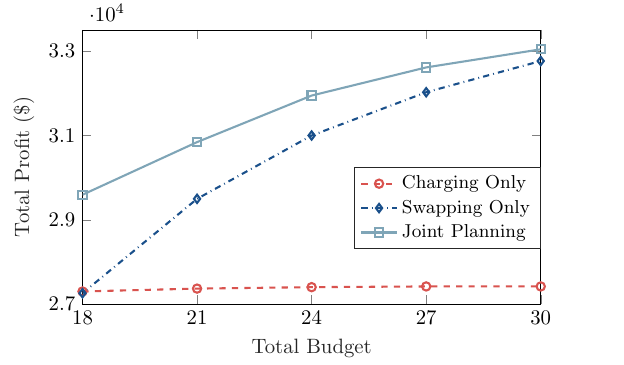}%
    \label{fig:cost ratio 3}}
    \hspace{-0.7cm}
    \subfloat[$\gamma_s/\gamma_c = 5$]{\includegraphics[width=0.35\textwidth]{figure/fig_total_profit.pdf}%
    \label{fig:cost ratio 5}}
    \hspace{-0.7cm}
    \subfloat[$\gamma_s/\gamma_c = 7$]{\includegraphics[width=0.35\textwidth]{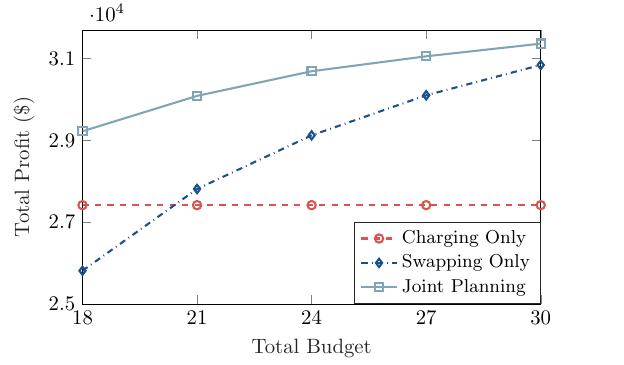}%
    \label{fig:cost ratio 7}}
    \caption{Performance under different cost ratios between battery swapping stations and charging stations.}
    \label{fig:sensitivity cost ratio}
\end{figure*}
\begin{figure*}[!th]
    \centering
    \subfloat[$\gamma_p = 15$]{\includegraphics[width=0.35\textwidth]{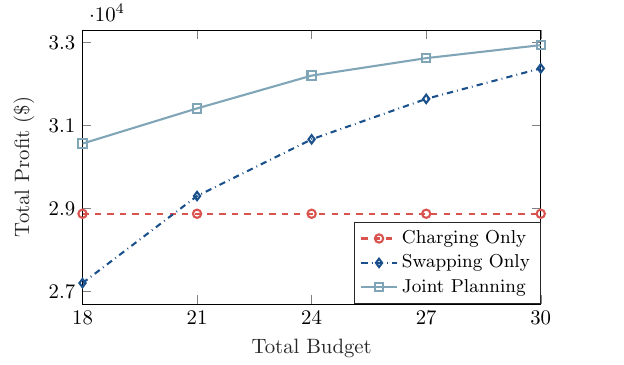}%
    \label{fig:gammap 15}}
    \hspace{-0.7cm}
    \subfloat[$\gamma_p = 20$]{\includegraphics[width=0.35\textwidth]{figure/fig_total_profit.pdf}%
    \label{fig:gammap 20}}
    \hspace{-0.7cm}
    \subfloat[$\gamma_p = 25$]{\includegraphics[width=0.35\textwidth]{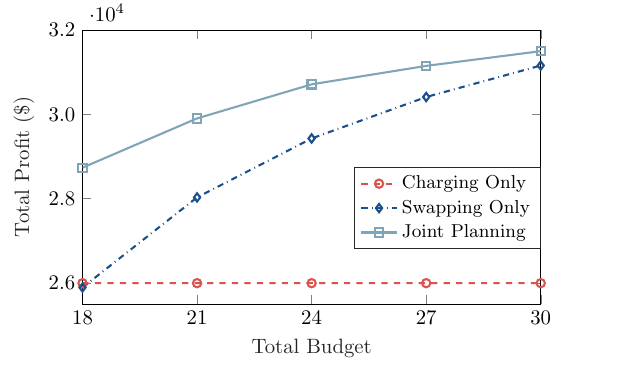}%
    \label{fig:gammap 25}}
    \caption{Performance under different values of the penalty parameter $\gamma_p$.}
    \label{fig:sensitivity gammap}
\end{figure*}
\begin{figure*}[!th]
    \centering
    \subfloat[$\gamma_e = 20$]{\includegraphics[width=0.35\textwidth]{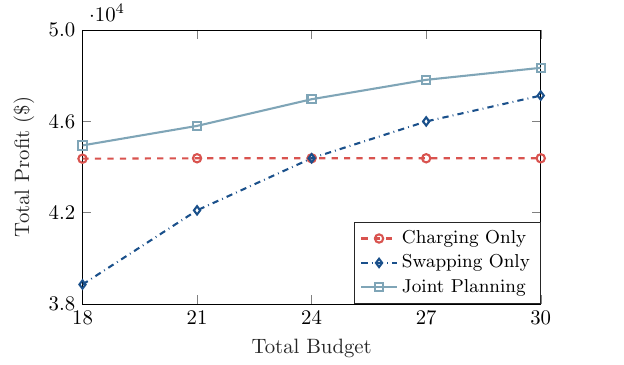}%
    \label{fig:gammae 20}}
    \hspace{-0.7cm}
    \subfloat[$\gamma_e = 25$]{\includegraphics[width=0.35\textwidth]{figure/fig_total_profit.pdf}%
    \label{fig:gammae 25}}
    \hspace{-0.7cm}
    \subfloat[$\gamma_e = 30$]{\includegraphics[width=0.35\textwidth]{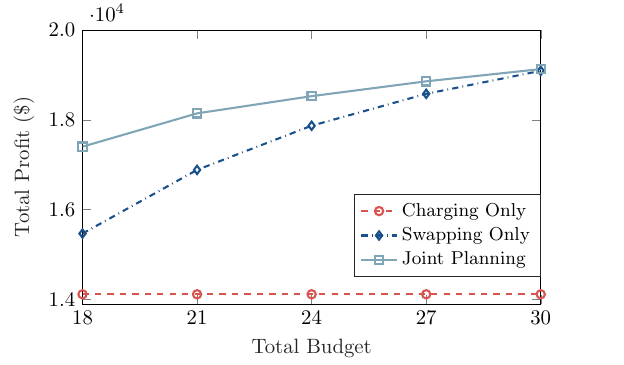}%
    \label{fig:gammae 30}}
    \caption{Performance under different values of the vehicle cost $\gamma_e$.}
    \label{fig:sensitivity gammae}
\end{figure*}
\begin{figure*}[!th]
    \centering
    \subfloat[$\gamma = 0.85$]{\includegraphics[width=0.35\textwidth]{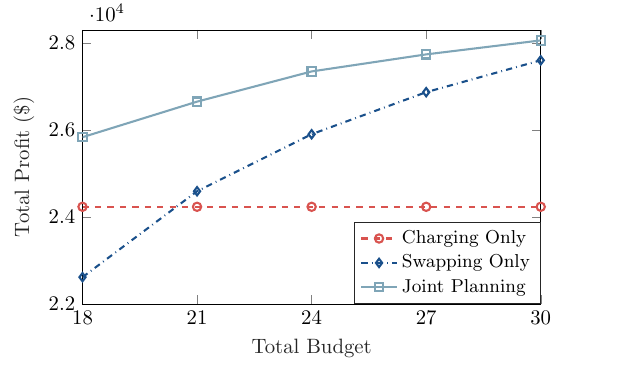}%
    \label{fig:discount 085}}
    \hspace{-0.7cm}
    \subfloat[$\gamma = 0.90$]{\includegraphics[width=0.35\textwidth]{figure/fig_total_profit.pdf}%
    \label{fig:discount 090}}
    \hspace{-0.7cm}
    \subfloat[$\gamma = 0.95$]{\includegraphics[width=0.35\textwidth]{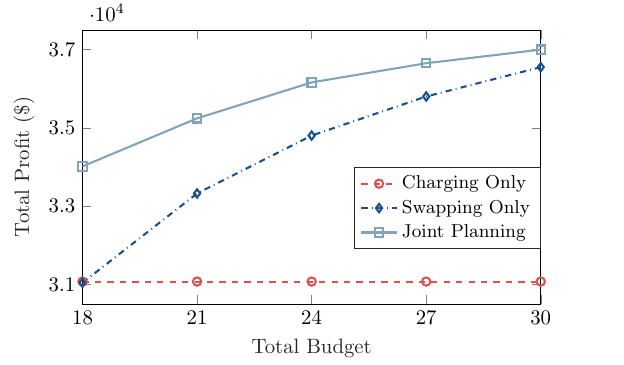}%
    \label{fig:discount 095}}
    \caption{Performance under different values of the discount factor $\gamma$.}
    \label{fig:sensitivity discount}
\end{figure*}
\begin{figure*}[!th]
    \centering
    \subfloat[$K = 30$]{\includegraphics[width=0.35\textwidth]{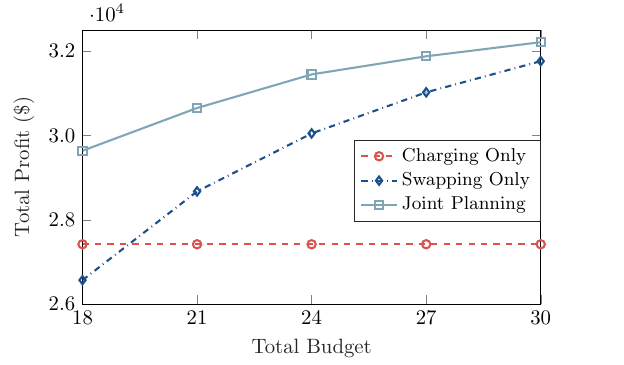}%
    \label{fig:K 30}}
    \hspace{-0.7cm}
    \subfloat[$K = 40$]{\includegraphics[width=0.35\textwidth]{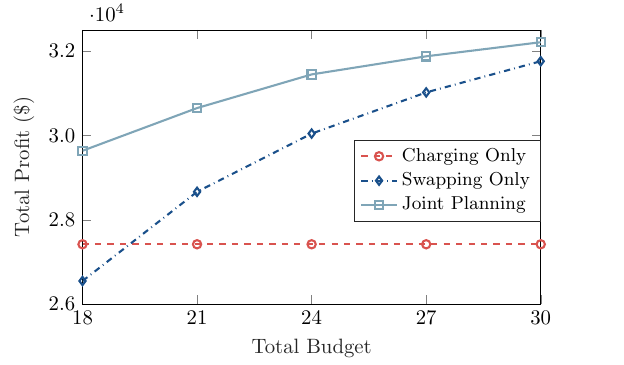}%
    \label{fig:K 40}}
    \hspace{-0.7cm}
    \subfloat[$K = 50$]{\includegraphics[width=0.35\textwidth]{figure/fig_total_profit.pdf}%
    \label{fig:K 50}}
    \caption{Performance under different values of the queue capacity $K$.}
    \label{fig:sensitivity K}
\end{figure*}
\begin{figure*}[!th]
    \centering
    \subfloat[$C = 3$]{\includegraphics[width=0.35\textwidth]{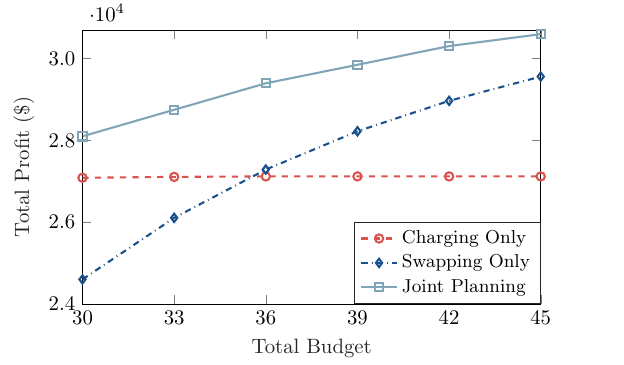}%
    \label{fig:C 3}}
    \hspace{-0.7cm}
    \subfloat[$C = 5$]{\includegraphics[width=0.35\textwidth]{figure/fig_total_profit}%
    \label{fig:C 5}}
    \hspace{-0.7cm}
    \subfloat[$C = 7$]{\includegraphics[width=0.35\textwidth]{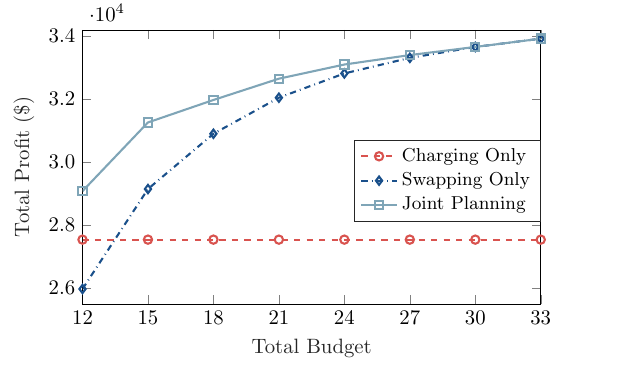}%
    \label{fig:C 7}}
    \caption{Performance under different infrastructure configurations.}
    \label{fig:sensitivity C}
\end{figure*}

\section{Results of Model Extensions}
\label{appendix:extension}
\begin{figure*}[!ht]
    \centering
    \subfloat[``First-Order'' Approximation]{\includegraphics[width=0.35\textwidth]{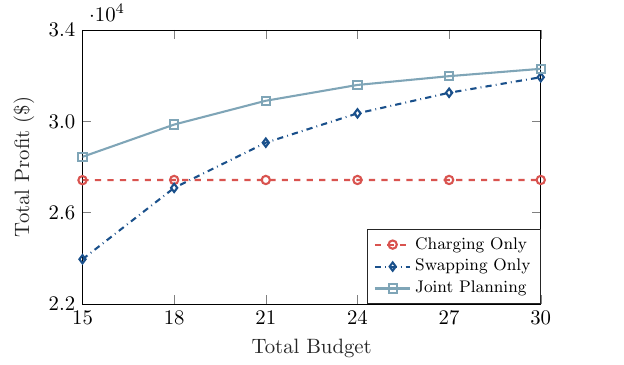}%
    \label{fig:blocking K 10 1}}
    \hspace{-0.7cm}
    \subfloat[``Second-Order'' Approximation]{\includegraphics[width=0.35\textwidth]{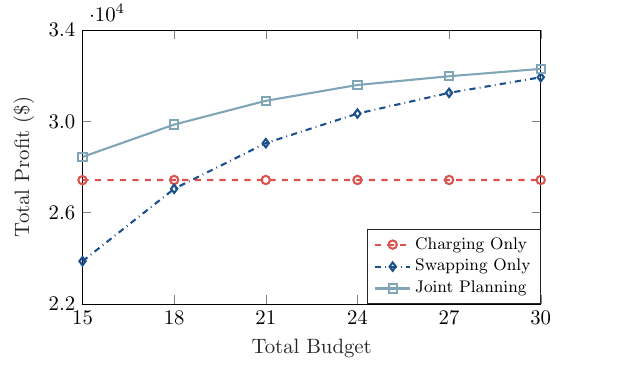}%
    \label{fig:blocking K 10 2}}
    \hspace{-0.7cm}
    \subfloat[Simplified Model]{\includegraphics[width=0.35\textwidth]{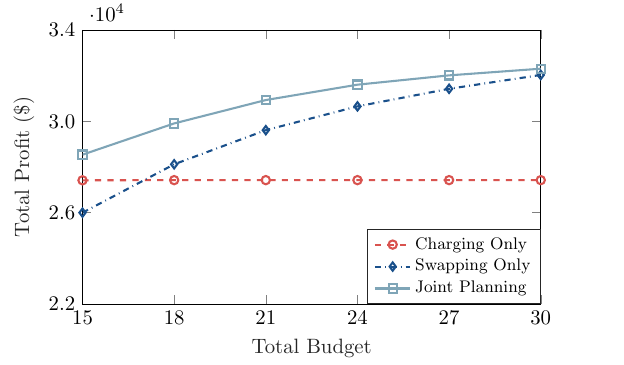}%
    \label{fig:blocking K 10 0}}
    \caption{Performance of the simplified and the extended models with a queue capacity of 10.}
    \label{fig:blocking K 10}
\end{figure*}
\begin{figure*}[!ht]
    \centering
    \subfloat[``First-Order'' Approximation]{\includegraphics[width=0.35\textwidth]{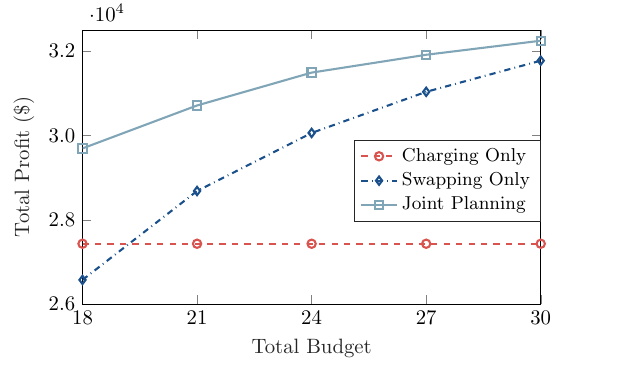}%
    \label{fig:blocking K 30 1}}
    \hspace{-0.7cm}
    \subfloat[``Second-Order'' Approximation]{\includegraphics[width=0.35\textwidth]{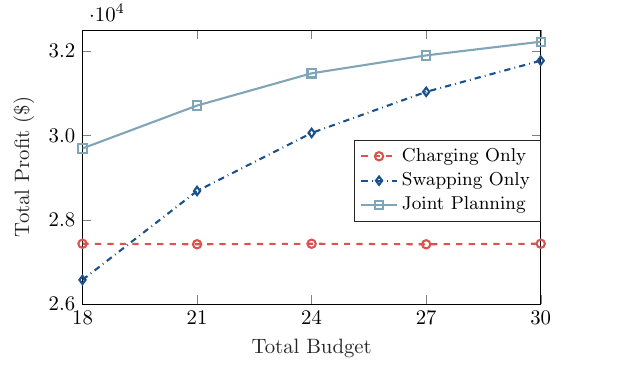}%
    \label{fig:blocking K 30 2}}
    \hspace{-0.7cm}
    \subfloat[Simplified Model]{\includegraphics[width=0.35\textwidth]{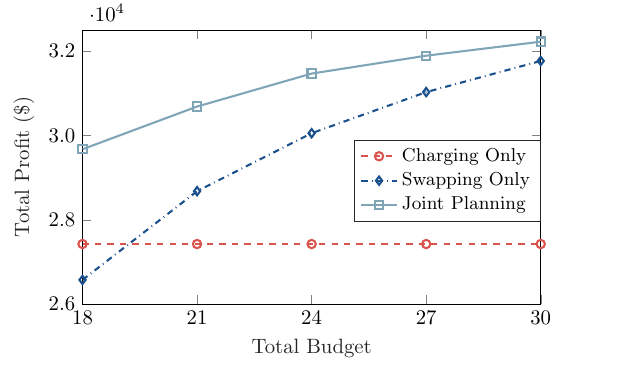}%
    \label{fig:blocking K 30 0}}
    \caption{Performance of the simplified and the extended models with a queue capacity of 30.}
    \label{fig:blocking K 30}
\end{figure*}

\begin{figure}[!th]
    \centering
    \includegraphics[width=0.42\textwidth]{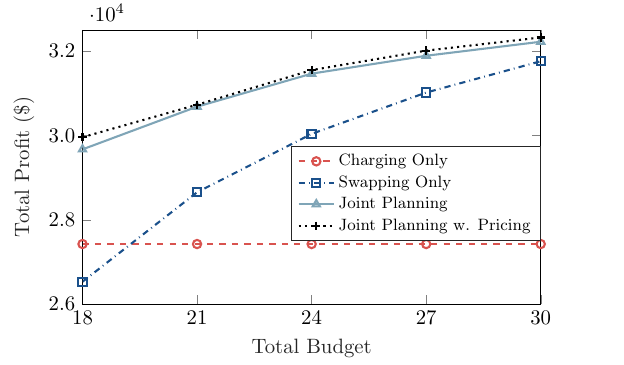}
    \caption{Comparison between the performance with and without differentiated pricing.}
    \label{fig:extension pricing}
\end{figure}

\end{document}